\documentclass[10pt,aps,pra,twocolumn,showpacs,showkeys,superscriptaddress,longbibliography]{revtex4-2}

\usepackage{times}
\usepackage{mathtools,amssymb,accents,bm}
\usepackage{siunitx}
\usepackage{hyperref}

\newlength{\dhatheight}
\newcommand{\hhatTS}[1]{%
    \settoheight{\dhatheight}{\ensuremath{\hat{#1}}}%
    \addtolength{\dhatheight}{-0.2ex}%
    \hat{\vphantom{\rule{1pt}{\dhatheight}}%
    \smash{\hat{#1}}}}    
\newcommand{\hhatS}[1]{%
    \settoheight{\dhatheight}{\ensuremath{\scriptstyle{\hat{#1}}}}%
    \addtolength{\dhatheight}{-0.15ex}%
    \hat{\vphantom{\rule{1pt}{\dhatheight}}%
    \smash{\hat{#1}}}}
\newcommand{\hhatSS}[1]{%
    \settoheight{\dhatheight}{\ensuremath{\scriptscriptstyle{\hat{#1}}}}%
    \addtolength{\dhatheight}{-0.075ex}%
    \hat{\vphantom{\rule{1pt}{\dhatheight}}%
    \smash{\hat{#1}}}}
\newcommand{\hhat}[1]{\mathchoice{\hhatTS{#1}}{\hhatTS{#1}}{\hhatS{#1}}{\hhatSS{#1}}}

\begin{document}

\title{Accelerated Newton-Raphson GRAPE methods for optimal control}

\author{David L. Goodwin}
\email[]{david.goodwin@partner.kit.edu}
\affiliation{Chemistry Research Laboratory, University of Oxford, Mansfield Road, Oxford OX1 3TA, UK}
\affiliation{Institute for Biological Interfaces 4 -- magnetic resonance, Karlsruhe Institute for Technology (KIT), Karlsruhe, Germany}

\author{Mads Sloth Vinding}
\email[]{msv@cfin.au.dk}
\affiliation{Center of Functionally Integrative Neuroscience (CFIN), Department of Clinical Medicine, Faculty of Health, Aarhus University, Denmark}

\date{\today}

\begin{abstract}
A Hessian based optimal control method is presented in Liouville space to mitigate previously undesirable polynomial scaling of computation time. This new method, an improvement to the state-of-the-art Newton-Raphson GRAPE method, is derived with respect to two exact time-propagator derivative techniques: auxiliary matrix and ESCALADE methods. We observed that compared to the best current implementation of Newton-Raphson GRAPE method, for an ensemble of 2-level systems, with realistic conditions, the new auxiliary matrix and ESCALADE Hessians can be 4-200 and 70-600 times faster, respectively.
\end{abstract}




\maketitle

\section{Introduction}

The problem of transferring the state of a dynamical system to a desired target state while minimising the remaining distance and costs is often solved with optimal control theory \cite{Glaser2015,Acin2018}. Applications include quantum sensing \cite{Saywell2018,SaywellThesis,Saywell2020,Saywell2021}, quantum computing \cite{Zhang2012,Waldherr2014,Dolde2014}, and nuclear magnetic resonance (\textsc{nmr}) spectroscopy \cite{Khaneja2001,Skinner2003,Khaneja2005,Tosner2006} and imaging (\textsc{mri}) \cite{Vinding2017,Reeth2019}.

A number of different approaches to optimal control (OC) has lead to the development of different methods: Lagrangian methods \cite{Somloi1993,Zhu1998,Maday2003,Eitan2011}; minimal-time OC \cite{Khaneja2001}; gradient ascent pulse engineering (\textsc{grape}) \cite{Khaneja2005}; sophisticated gradient-free searches \cite{Doria2011,Rach2015}; Krylov-Newton methods \cite{Ciaramella2015}; OC with a basis of analytic controls \cite{Machnes2018,Lucarelli2018}; a tensor product approach for large quantum systems \cite{QuinonesValles2019}. The method outlined in this work is based on a piecewise-constant control pulse approximation \cite{Conolly1986,Mao1986,Viola1998,Viola1999} of \textsc{grape} \cite{Khaneja2005,deFouquieres2011,Goodwin2016} using a gradient following numerical optimisation.

Although finding an optimal solution to the problem of controlling a single 2-level system from a defined initial system state to a desired target state, a \emph{state-to-state} problem, is considered straight-forward and computationally inexpensive with modern methods and computing power; an OC problem can become numerically and computationally arduous \cite{Pechen2011}, particularly for applications that account for practical hardware configurations and limitations \cite{Skinner2003,Kobzar2004,Kobzar2008}. Additionally, the computational expense can increase dramatically when optimising over an ensemble of systems, such as the case in solid-state nuclear magnetic resonance \cite{Tosner2006,GoodwinThesis,Tosner2021}, where the ensemble also includes crystalline orientations of a powder average \cite{Kuprov2016}. The particular application of OC of interest to the authors is that of a neural network based method for \textsc{mri} \cite{Vinding2019,Vinding2021a,Vinding2021b} and a method of morphic OC \cite{Goodwin2020,Haller2022}, requiring hundreds of thousands of optimised pulse shapes to form their optimal libraries.

Keeping the example of \textsc{mri}, the utility of OC is highlighted when considering the legal constraints of power deposition safeguards \cite{Dempsey2002} and the obvious financial rewards of reducing the time a patient stays inside the \textsc{mri} machine.

Modern techniques can mitigate protracted numerical convergence with a quadratically convergent optimisation method, requiring the calculation of a Hessian matrix, giving a large saving in the number of required serial optimisation iterations \cite{Goodwin2016,GoodwinThesis,Foroozandeh2021}. However, it is known within the community that calculation of the Hessian matrix does not scale well to a control problem with a large number of controllable amplitudes \cite{KalliesThesis} (also shown in Fig.~\ref{FIG_auxmat_timing}).

This communication presents a jump in computational efficiency of multiple-orders of magnitude with a rework of the original method, devised to calculate the exact Hessian matrix \cite{Goodwin2015,Goodwin2016,GoodwinThesis}, and the recently published exact, matrix-free, method of efficient spin control using analytical Lie algebraic derivatives (\textsc{escalade}) \cite{Foroozandeh2021}. The manuscript will present the mathematical formulation of these exact Newton-Raphson \textsc{grape} methods in the irreducible spherical-tensor basis of a Liouville space and show a newly devised method to calculate the Hessian matrix with $\mathcal{O}(N)$ scaling, compared to the previous $\mathcal{O}(N^2)$ scaling (Fig.~\ref{FIG_auxmat_timing}). Results show the comparative speedup of this method to the original in the context of state-to-state \textsc{mri} problems, requiring a generalisation of the \textsc{escalade} method with calculations of optimal $\mathrm{z}$-controls.

\section{Exact Newton-Raphson optimal control}

\subsection{Optimal control in Liouville space}

A quantum system can be described by a density operator, $\hat{\rho}(t)$; a time-dependent system state. The evolution of this state is dictated by the Liouville-von Neumann equation,
\begin{equation}
\frac{\partial}{\partial t} \hat{\rho}(t)=-i\big[\hat{H}(t),\hat{\rho}(t)\big]
\label{eq:LvN}
\end{equation}
where $\hat{H}(t)$ is a time-dependent Hamiltonian, an operator in a Hilbert space with a spectrum of allowed energy levels. The usual factor of $\hbar$ is dropped here, resulting in the eigenspectrum of $\hat{H}$ expressed in angular frequency units. The methods presented in this manuscript are particular to a Liouville space, also named the adjoint representation\cite{deFouquieres2011,Goodwin2016,GoodwinThesis}. A system state in a Liouville space is represented by a vector, $|\hat{\rho}\rangle$, obtained by stacking columns of the density operator $\hat{\rho}$ with Eq.~(\ref{eq:LvN}) becoming 
\begin{align}
& \frac{\partial}{\partial t} \big|\hat{\rho}(t)\big\rangle=-i\hhat{L}(t)\big|\hat{\rho}(t)\big\rangle, && \hhat{L}(t)\triangleq\openone\!\otimes\!\hat{H}(t)-\hat{H}(t)^\intercal\!\otimes\!\openone& 
\label{eq:liouville}
\end{align}
where the identity matrix, $\openone$, and Hamiltonian have the same dimension.

The form of a bilinear control problem is to split that which is controllable, \emph{the control}, from that which is not, \emph{the drift}. The Liouvillian for a control problem with $\mathrm{x}$-, $\mathrm{y}$-, and $\mathrm{z}$-controls on a 2-level system can be written as 
\begin{equation}
\hhat{L}(t)= \underbrace{\omega\hhat{L}_\mathrm{z}}_\text{drift} + \underbrace{c_\mathrm{x}(t)\hhat{L}_\mathrm{x}+c_\mathrm{y}(t)\hhat{L}_\mathrm{y}+c_\mathrm{z}(t)\hhat{L}_\mathrm{z}}_\text{control}\label{eq:HamTot}
\end{equation}
where the angular frequency $\omega$ is the time-independent resonant frequency offset, $c_\mathrm{x,y,z}(t)$ are time-dependent control amplitudes, and $\hhat{L}_\mathrm{x,y,z}$ are Pauli matrices of a Liouville space.

The OC method of \textsc{grape} \cite{Khaneja2005} uses piecewise constant control pulses, where control pulses are constant over a small time interval, $\Delta t$ \cite{Viola1998,Viola1999}:
\begin{align}
& c_{k}^{}(t)\to \begin{bmatrix} c_{k,1}^{} & c_{k,2}^{} & \cdots & c_{k,N}^{} \end{bmatrix}, && k\in\{\mathrm{x}, \mathrm{y}, \mathrm{z}\} & \label{eq:piecewise}
\end{align}
using the notation $c_{k,n}\equiv c_{k}(t_n)$ for convenience, and $t_N=N\Delta t$. This discrete formulation allows numerical solution to Eq.~(\ref{eq:liouville}), given an initial state of the system $|\hat{\rho}_0\rangle$, through time-ordered propagation
\begin{equation}
 |\hat{\rho}_N\rangle = \hhat{P}_{\!N}\hhat{P}_{\!N-1}\cdots\hhat{P}_{\!2}\hhat{P}_{\!1} |\hat{\rho}_0\rangle\label{eq:Uprop_liouv}
\end{equation}
where $\hhat{P}_{\!\!n}$ are time-propagators of an isolated time-slice and are defined through the exponential map
\begin{align}
&|\hat{\rho}_{n}\rangle = \hhat{P}_{\!\!n} |\hat{\rho}_{n-1}\rangle, && \hhat{P}_{\!\!n} = \mathrm{e}^{-i \hhat{L}_n \Delta t} &\label{eq:expmap}
\end{align}
The matrix exponential of Eq.~(\ref{eq:expmap}) usually calculated with the Pad\'{e} approximant, Taylor series, or Krylov propagation \cite{GoodwinThesis}.

An additional method to calculate the time-propagators of a 2-level system is by explicitly calculating the elements of the matrix. In a spherical-tensor basis, this matrix is the Wigner-matrix \cite{Siminovitch1997}:
\begin{equation}
\hhat{P}_{\!\!n}=\begin{bmatrix} \alpha^2 & \sqrt{2}\alpha\beta & \beta^2 \\ -\sqrt{2}\alpha\beta^\ast & \alpha\alpha^\ast-\beta\beta^\ast & \sqrt{2}\alpha^\ast\beta \\ {\beta^\ast}^2 & -\sqrt{2}\alpha^\ast \beta^\ast & {\alpha^\ast}^2 \end{bmatrix}\label{EQ_Wigner}
\end{equation}
which is formulated in terms of the complex elements 
\begin{align}
& \alpha = \cos{\phi}+i\frac{z}{r} \sin{\phi}, && \beta = \frac{y}{r} \sin{\phi} + i\frac{x}{r} \sin{\phi}, &\label{EQ:alphabeta}
\end{align}
where the shorthands $x=c_{\mathrm{x},n}$, $y=c_{\mathrm{y},n}$, and $z=c_{\mathrm{z},n}+\omega$ are used, $\phi= \frac{1}{2}r\Delta t$ is a polar angle of rotation, and $r=\sqrt{x^2+y^2+z^2}$ is the polar radius.

As a notational convenience for what follows, the following \emph{effective propagators} are defined as the effect of the pulse between time slices $m$ and $n$:
\begin{align}
& \mathbf{U}_{m}^{n} \triangleq \hhat{P}_{\!\!n} \hhat{P}_{\!\!n-1} \cdots \hhat{P}_{\!\!m+1} \hhat{P}_{\!\!m}, && \forall (1\leqslant m\leqslant n\leqslant N) & \label{eq:eff_prop}
\end{align}
where a backward (time-reversed) propagation can be denoted by $\mathbf{U}_{n}^{m}={\mathbf{U}_{m}^{n}}^\dagger$.

Optimal control requires a metric to optimise; the \emph{fidelity} \cite{Glaser1998,Khaneja2005}, $\mathcal{F}$, a measure how well the pulses perform a desired control task. The task of the OC problem is to find a set of control amplitudes, $c_\mathrm{x,y,z}(t)$, that maximise the fidelity e.g. the real part of an inner product:
\begin{equation}
\max_{c_\mathrm{x,y,z}(t)}\!\big(\mathcal{F}\big) = \max_{c_\mathrm{x,y,z}(t)}\!\Big(\mathrm{Re}\big\langle \hat{\sigma} \big| \mathbf{U}_{1}^{N}|\hat{\rho}_0 \big\rangle\Big)\label{eq:fidelity_vector}
\end{equation}
This form of the fidelity metric is defined in terms of state-to-state problems, where a system is in a defined initial state, $| \hat{\rho}_0 \rangle$, and the control task is to take this state to a desired one, $| \hat{\sigma} \rangle$. With the notation introduced in Eq.~(\ref{eq:eff_prop}), $\mathbf{U}_1^N$ is interpreted as the effective propagator over the shaped pulse. In addition to Eq.~(\ref{eq:Uprop_liouv}), the system is propagated backwards from the desired target state,
\begin{equation}
 |\hat{\chi}_n\rangle = \hhat{P}_{\!\!n}^\dagger\hhat{P}_{\!\!n+1}^\dagger\cdots\hhat{P}_{\!N-1}^\dagger\hhat{P}_{\!N}^\dagger |\hat{\sigma}\rangle\label{eq:Upropadj_liouv}
\end{equation}
which is termed the \emph{adjoint state} of the control problem \cite{GoodwinThesis}.

\textsc{Grape} is a gradient following numerical optimisation method and requires derivatives of the fidelity with respect to the controls. In turn, this requires the directional propagator derivatives, $D_{\!k,n}^a$; with the subscripts $k$ denoting the derivative in the direction of $\hhat{L}_k$ and $n$ denoting the derivative operating on the time-propagator $\hhat{P}_{\!\!n}$. For each time-slice, $n$, and for each control direction $\hhat{L}_k\in\{\hhat{L}_\mathrm{x},\hhat{L}_\mathrm{y},\hhat{L}_\mathrm{z}\}$, the $a^\text{th}$ order derivative takes the form
\begin{equation}
\nabla_{}^{a}\!\mathcal{F}(c_{k,n}^{}) = \mathrm{Re}\underbrace{\big\langle\hat{\chi}_{n+1}}_{\mathclap{N \text{ times}}}\!\big|\overbrace{\!D_{\!k,n}^a\hat{\rho}_{n-1}\big\rangle}^{\mathclap{N \text{ times}}},\label{eq:gradF}
\end{equation}
where bra-ket notation explicitly shows vector structures i.e. $|D_{\!k,n}^a \hat{\rho}_{n-1}\rangle=D_{\!k,n}^a|\hat{\rho}_{n-1}\rangle$. The number of forward and backward propagations is indicated for each control channel to produce a gradient vector, $\nabla\!\mathcal{F}$, or the diagonal elements of a Hessian matrix, $\nabla^2\!\mathcal{F}$.

\begin{figure}
\centering{\includegraphics{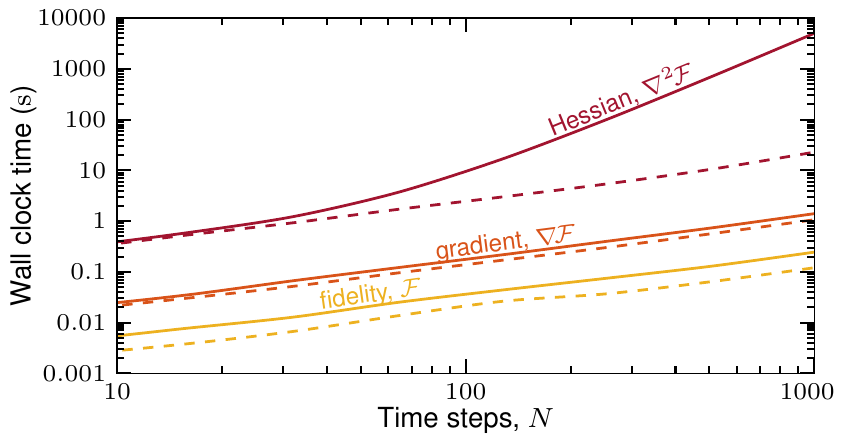}}
\caption{Average wall-clock time of fidelity, gradient, and Hessian for increasing $N$. Solid lines show the Newton-Raphson \textsc{grape} method \cite{Goodwin2016} and dashed lines show the proposed \emph{accelerated} Newton-Raphson \textsc{grape} method. Both methods use the auxiliary matrix method \cite{Goodwin2015} to calculate propagator derivatives.}
\label{FIG_auxmat_timing}
\end{figure}

This is sufficient for a fidelity gradient, scaling linearly with $N$, but a fidelity Hessian also requires mixed second order derivatives \cite{Goodwin2016}, where the off-diagonal Hessian elements are
\begin{equation}
\nabla_{}^{2}\!\mathcal{F}(c_{k,n}^{},c_{j,m}^{}) = \mathrm{Re}\underbrace{\!\big\langle\hat{\chi}_{n+1} {D_{\!k,n}^{}}}_{\mathclap{N \text{ times}}}\!\big|\overbrace{\!\mathbf{U}_{m+1}^{n-1}\big|D_{\!j,m}^{}\hat{\rho}_{m-1}\big\rangle}^{\mathclap{\frac{1}{2}N(N-1) \text{ times}}}\label{eq:HessF21}
\end{equation}
Clearly, the form of Eq.~(\ref{eq:HessF21}) has a central propagator that cannot be absorbed into the bra or ket because it depends on both $t_n$ and $t_m$, and therefore the computation scales polynomially with $\frac{1}{2}N(N-1)$ (the factor $\frac{1}{2}$ comes from the symmetric property of a Hessian \cite{GoodwinThesis}). This is known within the OC community \cite{KalliesThesis} and is highlighted in Fig.~\ref{FIG_auxmat_timing}. The linear plots of the fidelity and gradient, on log-log axes, show these calculations are efficient with increasing $N$, whereas the Hessian calculation time is not linear on these log-log axes. The subject of this manuscript is to mitigate this undesirable scaling, resulting in a linearly scaling Hessian calculation (dashed lines in Fig.~\ref{FIG_auxmat_timing}), after the following section outlines calculation of directional propagator derivatives.

\subsection{Directional propagator derivatives with auxiliary matrices}

As has been published previously \cite{Goodwin2015,Goodwin2016,GoodwinThesis}, exact propagator derivatives required by Eq.~(\ref{eq:gradF}) and (\ref{eq:HessF21}) can be calculated by exponentiating an auxiliary matrix: resulting in an upper triangular block matrix \cite{VanLoan1978} with a time-propagator, $\hhat{P}_{\!\!n}$, on the block diagonal, and with the directional derivative of that propagator, calculated in the direction of a control operator, in the upper triangular block. The first and second order propagator derivatives can be extracted from
\begin{equation}
\exp{\begin{bmatrix} \mathbf{A}_n & \mathbf{C}_j & \bm{0} \\ \bm{0} & \mathbf{A}_n & \mathbf{C}_k \\ \bm{0} & \bm{0} & \mathbf{A}_n\end{bmatrix}}=\begin{bmatrix} \hhat{P}_{\!\!n} &D_{\!j,n}^{} & \frac{1}{2}D_{\!jk,n}^2 \\ \bm{0} & \hhat{P}_{\!\!n} &D_{\!k,n}^{} \\ \bm{0} & \bm{0} & \hhat{P}_{\!\!n}  \end{bmatrix}\label{eq:auxmat}
\end{equation}
where the block matrix is formed from $\mathbf{A}_n=-i\hhat{L}_n\Delta t$, as in Eq.~(\ref{eq:expmap}), and $\mathbf{C}_k=-i\hhat{L}_k\Delta t$ are functions of the control operators with $j,k\in\{\mathrm{x},\mathrm{y},\mathrm{z}\}$.

\subsection{Directional propagator derivatives with \textsc{escalade}}

Whereas Eq.~(\ref{eq:auxmat}) is exact, the matrix exponential is expensive. Foroozandeh and Singh recently devised a method that is free from this expensive matrix exponential calculation, in the OC method of Efficient Spin Control using Analytical Lie Algebraic Derivatives (\textsc{escalade}) \cite{Foroozandeh2021}. This method can also be extended with interaction propagator splitting \cite{Goodwin2022}.

The efficient method of \textsc{escalade} to calculate the directional derivatives in Eqs.~(\ref{eq:gradF}) and (\ref{eq:HessF21}), in this single spin model, is to construct a matrix with rows containing all elements needed to construct the propagator derivatives    
\begin{equation}
\begin{bmatrix} \mathbf{\Theta}_{\!\mathrm{x}}^{} \\ \mathbf{\Theta}_{\!\mathrm{y}}^{} \\ \mathbf{\Theta}_{\!\mathrm{z}}^{} \end{bmatrix} = \mathrm{vec}[\openone] + \frac{\sin^2{\phi}}{\phi r}\mathrm{vec}[\mathbf{S}] + \frac{2\phi - \sin{2\phi}}{2\phi r^2}\mathrm{vec}[\mathbf{S}^2]\label{EQ_Dmatrix}
\end{equation}
 where $\mathrm{vec}[\openone]$ is a vectorised identity matrix, with a vectorisation operation on a matrix $\mathbf{A}$ such that $\mathbf{A}=\mathrm{vec}^{-1}[\mathrm{vec}[\mathbf{A}]]$. The skew-symmetric matrix $\mathbf{S}$ is
 \begin{equation}
\mathbf{S}= \begin{bmatrix} 0 & z & -y \\ -z & 0 & x \\ y & -x & 0  \end{bmatrix}\label{eq:Smats}
\end{equation}
and the symmetric matrix $\mathbf{S}^2$ is calculated algebraically. In turn, the directional propagator derivatives can be written as
\begin{gather}
D_{\!k,n}^{}= \hhat{P}_{\!\!n}\Big[\mathrm{vec}^{-1}\big[\mathbf{\Sigma}\mathbf{\Theta}_{\!k}^{}\big]\Big]
\label{EQ_dPdc}\\
\begin{aligned}
D_{\!jk,n}^{2}= & \hhat{P}_{\!\!n}\Big[\mathrm{vec}^{-1}\big[\mathbf{\Sigma}\mathbf{\Theta}_{\!jk}^{}\big]\Big]+\\ & \hhat{P}_{\!\!n}\Big[\mathrm{vec}^{-1}\big[\mathbf{\Sigma}\mathbf{\Theta}_{\!j}^{}\big]\Big]\Big[\mathrm{vec}^{-1}\big[\mathbf{\Sigma}\mathbf{\Theta}_{\!k}^{}\big]\Big]
\end{aligned}
\label{EQ_d2Pdc2}
\end{gather}
where $\mathbf{\Sigma}= \begin{bmatrix}\mathrm{vec}[\mathbf{C}_\mathrm{x}] & \mathrm{vec}[\mathbf{C}_\mathrm{y}] & \mathrm{vec}[\mathbf{C}_\mathrm{z}]\end{bmatrix}$ is a 3-column matrix to be multiplied with the 3-row matrix of $\mathbf{\Theta}_{\!k}^{}$. Using a similar notation to Eq.~(\ref{EQ_Dmatrix}), the $\mathbf{\Theta}_{\!jk}^{}$ matrices required by the second order derivatives in Eq.~(\ref{EQ_d2Pdc2}) are
\begin{align}
\begin{bmatrix} \mathbf{\Theta}_{\!j\mathrm{x}}^{} \\ \mathbf{\Theta}_{\!j\mathrm{y}}^{} \\ \mathbf{\Theta}_{\!j\mathrm{z}}^{} \end{bmatrix} =   &  \frac{\cos(2\phi) -1}{r^2} \mathrm{vec}\bigg[\mathbf{S}\frac{\mathrm{d}\mathbf{S}}{\mathrm{d}j}+\frac{\mathrm{d}\mathbf{S}}{\mathrm{d}j}\mathbf{S}\bigg] \nonumber\\
              & + \frac{2\phi -\sin(2\phi)}{r^3} \mathrm{vec}\bigg[\frac{\mathrm{d}\mathbf{S}}{\mathrm{d}j}\bigg] \nonumber\\
              & + \frac{2k(1-\cos(2\phi)-\phi\sin(2\phi))}{r^3} \mathrm{vec}[\mathbf{S}] \nonumber\\
              & + \frac{3k\sin(2\phi)-2k\phi(2+\cos(2\phi))}{r^4} \mathrm{vec}[\mathbf{S}^2]
\end{align}
where three such equations are required for $j\in\{\mathrm{x},\mathrm{y},\mathrm{z}\}$ and the derivatives of $\mathbf{S}$ are derived algebraically from Eq.~(\ref{eq:Smats}).

Since the directional derivatives of Eq.~(\ref{EQ_dPdc}) and Eq.~(\ref{EQ_d2Pdc2}), together with time-propagators of Eq.~(\ref{EQ_Wigner}), do not involve any matrix operations, other than a few trivial multiplications, \textsc{escalade} offers substantial computational gains relative to the \textsc{auxmat} method in Eq.~(\ref{eq:auxmat}).

\section{Accelerated Newton-Raphson optimal control}

Moving away from a chosen calculation method of directional propagator derivatives, the remaining bottleneck of the Newton-Raphson methods presented above is the off-diagonal Hessian elements; the mixed derivatives of Eq.~(\ref{eq:HessF21}).

As a starting point, \textsc{escalade} includes an additional efficiency, which also applies to the \textsc{auxmat} method, where the central effective propagator in Eq.~(\ref{eq:HessF21}) can be split to $\mathbf{U}_{m+1}^{n-1}=\mathbf{U}_{m+1}^{n-1}\mathbf{U}_1^{m}\mathbf{U}_m^{1}=\mathbf{U}_{1}^{n-1}\mathbf{U}_m^{1}$, Eq.~(\ref{eq:HessF21}) becomes
\begin{equation}
\nabla_{}^{2}\!\mathcal{F}(c_{k,n}^{},c_{j,m}^{}) = \mathrm{Re}\big\langle\hat{\chi}_{n+1} D_{\!k,n}^{}\big|\mathbf{U}_{1}^{n-1}\mathbf{U}_m^{1}\big|D_{\!j,m}^{}\hat{\rho}_{m-1}\big\rangle\label{EQ_esc_prop_trick}
\end{equation}
An interpretation of these two central effective propagators is: the right hand side, $\langle\hat{\chi}_{n+1} D_{\!k,n}^{}|$, is multiplied by the effective propagator $\mathbf{U}_{1}^{n-1}$, presenting the directional derivative of $t_n$ to be evaluated at $t_0$; the left hand side, $|D_{\!j,m}^{}\hat{\rho}_{m-1}\rangle$, is multiplied by the time-reversed effective propagator $\mathbf{U}_m^{1}$, presenting the directional derivative of $t_m$ to also be evaluated at $t_0$.

To outline how this can be overcome, a representation of a \emph{trajectory} is introduced:
\begin{equation}
\big[\bm{\rho}\big]_m^n\triangleq\begin{bmatrix}|\hat{\rho}_n\rangle & |\hat{\rho}_{n-1}\rangle & \cdots & |\hat{\rho}_{m+1}\rangle & |\hat{\rho}_m\rangle\end{bmatrix}
\end{equation}
which is an array of column vectors and the whole trajectory from the dynamics in Eq.~(\ref{eq:HamTot}) is contained in $[\bm{\rho}]_0^N$. Trajectory analysis is useful for visualising pulse dynamics \cite{Kuprov2013}. Taking the concept of a single matrix containing all trajectory information, a \emph{directional derivative trajectory}, evaluated at $t_0$ as with the right side of Eq.~(\ref{EQ_esc_prop_trick}), can be defined as
\begin{equation}
\big[\partial_j^{[0]}\!\bm{\rho}\big]_0^n\triangleq\begin{bmatrix}|\mathbf{U}_n^{1}D_{\!j,n}^{}\hat{\rho}_{n-1}\rangle & \cdots & |\mathbf{U}_2^{1}D_{\!j,2}^{}\hat{\rho}_{1}\rangle & |\mathbf{U}_1^{\dagger}D_{\!j,1}^{}\hat{\rho}_{0}\rangle\end{bmatrix}\label{EQ_derivarray}
\end{equation}
where the superscript $^{[0]}$ is used to indicate evaluation at $t_0$.

With the realisation that a \emph{directional derivative trajectory} in Eq.~(\ref{EQ_derivarray}) is a matrix in itself, $n$ Hessian elements can be calculated with one matrix-vector product with
\begin{equation}
\begin{bmatrix}\nabla_{}^{2}\!\mathcal{F}(c_{k,n}^{},c_{j,n}^{}) \\ \vdots \\ \nabla_{}^{2}\!\mathcal{F}(c_{k,n}^{},c_{j,2}^{})  \\ \nabla_{}^{2}\!\mathcal{F}(c_{k,n}^{},c_{j,1}^{}) \end{bmatrix}^{\!\mathrm{T}}\!\!\!\!= \mathrm{Re}\underbrace{\big\langle\hat{\chi}_{n+1} D_{\!k,n}^{}\mathbf{U}_{1}^{n-1}}_{\mathclap{N \text{ times backward}}}\!\big|\overbrace{\!\big[\partial_j^{[0]}\!\bm{\rho}\big]_0^n}^{\mathclap{N \text{ times forward}}}\label{EQ_accMethod}
\end{equation}

Given that the total effective propagator $\mathbf{U}_{1}^{N}$ can be calculated with forward propagation, $\mathbf{U}_{1}^{n-1}$ can be updated during subsequent backward propagation with $\mathbf{U}_{1}^{n-1}=\hhat{P}_{\!\!n-1}^\dagger\mathbf{U}_{1}^{n}$. The single vector-matrix product of Eq.~(\ref{EQ_accMethod}), per time-slice, is expected to be much more efficient than the $n-1$ vector-vector products, per time-slice, in Eq.~(\ref{eq:HessF21}).

\section{Acceleration of a broadband MRI example}

\begin{figure*}
\centering{\includegraphics{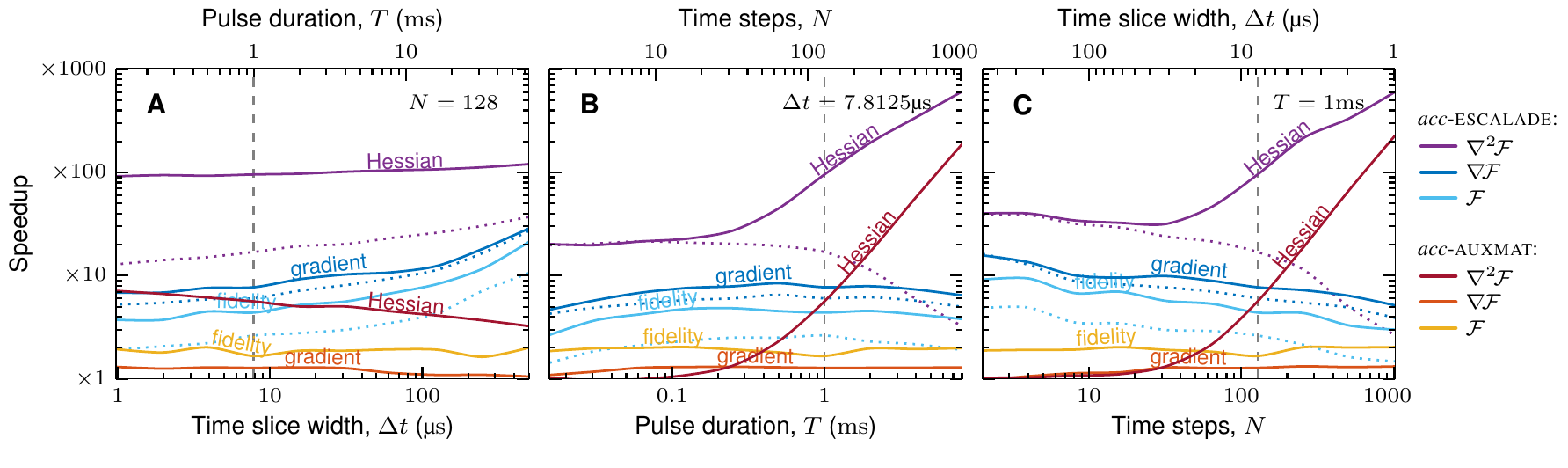}}
\caption{The average speedup of the accelerated \textsc{escalade} method (\emph{acc}-\textsc{escalade}) and the accelerated auxiliary matrix method (\emph{acc}-\textsc{auxmat}), compared to the auxiliary matrix method. The number of time-steps, $N$, the time-slice width, $\Delta t$, and the pulse duration, $T$, are parameters affecting computation time. In each of the three plots one parameter is kept constant: (A) $N=128$, (B) $\Delta t=7.8125$~\si{\micro\second}; and (C) $T=1$~\si{\milli\second}. The dotted lines in each plot show the speedup of \emph{acc}-\textsc{escalade} relative to \emph{acc}-\textsc{auxmat}.}
\label{FIG_speedup}
\end{figure*}

The authors choose to set the test of the method presented above in the context of \textsc{mri} because there is a very real need to have computationally fast OC methods to run \emph{on-the-fly} with a need for patient specific solutions. Often, OC-facilitated \textsc{mri} exploits simultaneous/parallel control systems e.g. eight radiofrequency channels ($\mathrm{x}$- and $\mathrm{y}$-controls) \cite{Vinding2017,vinding_local_2017} and three magnetic field gradients and/or numerous local shims ($\mathrm{z}$-controls) \cite{vinding_local_2017,stockmann_vivo_2018,Vinding2021b}.

The optimisation is set as robust for an ensemble of $101$ spin-$\frac{1}{2}$ systems, with offset bandwidth $\frac{\omega}{2\pi}\in[-\frac{1}{2},+\frac{1}{2}]~\si{\kilo\hertz}$, simulated using a single block-diagonal Liouvillian with each block operating on a single ensemble member \cite{Kuprov2016,Goodwin2020}.

The time-slice, $\Delta t$, is often fixed to the hardware digitalisation dwell-time: typically a few \si{\micro\second}. Modern \textsc{mri} systems handle $N>1000$, but the pulse duration, $T= N\Delta t$, is usually kept around a few \si{\milli\second} to avoid degradation of pulse performance due to transverse relaxation. As a set of speed tests, these three variables are incremented and set as a range of OC problems and are shown in Fig.~\ref{FIG_speedup}. The vertical dashed lines indicate a set of variables physically relevant to \textsc{mri}.

Figure~\ref{FIG_speedup} shows the average speedup of the accelerated auxiliary matrix method (\emph{acc}-\textsc{auxmat}) and the accelerated \textsc{escalade} method (\emph{acc}-\textsc{escalade}), both relative to the standard auxiliary matrix method. A separate fidelity, gradient, and Hessian calculation is performed for each method from a random control waveform, with $c_\mathrm{x}$, $c_\mathrm{y}$, and $c_\mathrm{z}$ controls in the order of $2\pi\times10^3~\si{\radian\per\second}$. The average is over 84 different random control pulses, run in parallel on 28 \textsc{cpu} cores.

Figure~\ref{FIG_speedup}\textbf{A} shows that the effect of increasing $T$ while also increasing $\Delta t$ is, approximately, a constant speedup of \emph{acc}-\textsc{auxmat} and \emph{acc}-\textsc{escalade} Hessian. This can be attributed to the effective propagator splitting of Eq.~(\ref{EQ_esc_prop_trick}). The small speedup of \emph{acc}-\textsc{auxmat} fidelity and gradient calculations is due to an algorithmic efficiency from coding in a way that lends itself to the accelerated Hessian methods. 

The main result of this manuscript is the speedup of Hessian calculations of the \emph{acc}-\textsc{auxmat} and \emph{acc}-\textsc{escalade} methods, with speedup increasing as $N$ increases, in Fig.~\ref{FIG_speedup}\textbf{B} and Fig.~\ref{FIG_speedup}\textbf{C}, to over $\times 100$ when a large $N$ is used. Furthermore, this trend does not appear to dampen, indicating further speedup when $N\gg 1000$.

From the evidence in Fig.~\ref{FIG_speedup}\textbf{B} and Fig.~\ref{FIG_speedup}\textbf{C}, the accelerated Hessian calculations of Eq.~(\ref{EQ_accMethod}) do indeed remove the $\mathcal{O}(N^2)$ scaling of the original Newton-Raphson method \cite{Goodwin2016}, reducing to a linear scaling $\mathcal{O}(N)$.

\section{Conclusion}

A new mathematical formulation of the Newton-Raphson \textsc{grape} method \cite{Goodwin2016} has been presented in Liouville space and applies to optimal control problems with unitary evolution. The \textsc{escalade} method \cite{Foroozandeh2021} recast the cumbersome problem of finding derivatives -- using trigonometric evaluations of vectorised arrays rather than matrix-matrix products and computationally expensive matrix exponentials -- and factorised the central propagator to avoid $\mathcal{O}(N^2)$ scaling, now reduced to $\mathcal{O}(N)$, in the computation of the Hessian. The \textsc{escalade} method is shown in a new light of Liouville space, with the additional derivation of $\mathrm{z}$-controls, which are important to \textsc{mri}. 

The main result is a new formulation of the Newton-Raphson \textsc{grape} method in Liouville space, being an optimisation method with quadratic convergence to an optimal solution, reducing the expensive polynomial scaling of the control problem to a linear scaling when increasing the number of piecewise-constant pulses in optimal solution. Speedup increases as the number of time-slices, $N$, increases: $\times 4$ to $\times 200$, for $N=100$ and $N=1000$ respectively. Furthermore, employing \textsc{escalade} within this new Hessian calculation method shows further speedup of $\times 70$ to $\times 600$, for $N=100$ and $N=1000$ respectively, when compared to the original Newton-Raphson \textsc{grape} method.

\begin{acknowledgments}
DLG thanks Pranav Singh and Mohammadali Foroozandeh for useful discussions on the \textsc{escalade} method, Burkhard Luy for his continued support, and Martin Goodwin for proofreading.

MSV would like to thank Villum Fonden, Eva og Henry Fraenkels Mindefond, Harboefonden, and Kong Christian Den Tiendes Fond.
\end{acknowledgments}


\begin{thebibliography}{52}%
\makeatletter
\providecommand \@ifxundefined [1]{%
 \@ifx{#1\undefined}
}%
\providecommand \@ifnum [1]{%
 \ifnum #1\expandafter \@firstoftwo
 \else \expandafter \@secondoftwo
 \fi
}%
\providecommand \@ifx [1]{%
 \ifx #1\expandafter \@firstoftwo
 \else \expandafter \@secondoftwo
 \fi
}%
\providecommand \natexlab [1]{#1}%
\providecommand \enquote  [1]{``#1''}%
\providecommand \bibnamefont  [1]{#1}%
\providecommand \bibfnamefont [1]{#1}%
\providecommand \citenamefont [1]{#1}%
\providecommand \href@noop [0]{\@secondoftwo}%
\providecommand \href [0]{\begingroup \@sanitize@url \@href}%
\providecommand \@href[1]{\@@startlink{#1}\@@href}%
\providecommand \@@href[1]{\endgroup#1\@@endlink}%
\providecommand \@sanitize@url [0]{\catcode `\\12\catcode `\$12\catcode
  `\&12\catcode `\#12\catcode `\^12\catcode `\_12\catcode `\%12\relax}%
\providecommand \@@startlink[1]{}%
\providecommand \@@endlink[0]{}%
\providecommand \url  [0]{\begingroup\@sanitize@url \@url }%
\providecommand \@url [1]{\endgroup\@href {#1}{\urlprefix }}%
\providecommand \urlprefix  [0]{URL }%
\providecommand \Eprint [0]{\href }%
\providecommand \doibase [0]{https://doi.org/}%
\providecommand \selectlanguage [0]{\@gobble}%
\providecommand \bibinfo  [0]{\@secondoftwo}%
\providecommand \bibfield  [0]{\@secondoftwo}%
\providecommand \translation [1]{[#1]}%
\providecommand \BibitemOpen [0]{}%
\providecommand \bibitemStop [0]{}%
\providecommand \bibitemNoStop [0]{.\EOS\space}%
\providecommand \EOS [0]{\spacefactor3000\relax}%
\providecommand \BibitemShut  [1]{\csname bibitem#1\endcsname}%
\let\auto@bib@innerbib\@empty
\bibitem [{\citenamefont {Glaser}\ \emph {et~al.}(2015)\citenamefont {Glaser},
  \citenamefont {Boscain}, \citenamefont {Calarco}, \citenamefont {Koch},
  \citenamefont {K{\"{o}}ckenberger}, \citenamefont {Kosloff}, \citenamefont
  {Kuprov}, \citenamefont {Luy}, \citenamefont {Schirmer}, \citenamefont
  {Schulte-Herbr{\"{u}}ggen}, \citenamefont {Sugny},\ and\ \citenamefont
  {Wilhelm}}]{Glaser2015}%
  \BibitemOpen
  \bibfield  {author} {\bibinfo {author} {\bibfnamefont {S.~J.}\ \bibnamefont
  {Glaser}}, \bibinfo {author} {\bibfnamefont {U.}~\bibnamefont {Boscain}},
  \bibinfo {author} {\bibfnamefont {T.}~\bibnamefont {Calarco}}, \bibinfo
  {author} {\bibfnamefont {C.~P.}\ \bibnamefont {Koch}}, \bibinfo {author}
  {\bibfnamefont {W.}~\bibnamefont {K{\"{o}}ckenberger}}, \bibinfo {author}
  {\bibfnamefont {R.}~\bibnamefont {Kosloff}}, \bibinfo {author} {\bibfnamefont
  {I.}~\bibnamefont {Kuprov}}, \bibinfo {author} {\bibfnamefont
  {B.}~\bibnamefont {Luy}}, \bibinfo {author} {\bibfnamefont {S.}~\bibnamefont
  {Schirmer}}, \bibinfo {author} {\bibfnamefont {T.}~\bibnamefont
  {Schulte-Herbr{\"{u}}ggen}}, \bibinfo {author} {\bibfnamefont
  {D.}~\bibnamefont {Sugny}},\ and\ \bibinfo {author} {\bibfnamefont {F.~K.}\
  \bibnamefont {Wilhelm}},\ }\bibfield  {title} {\bibinfo {title} {Training
  schr{\"{o}}dinger's cat: quantum optimal control},\ }\href
  {https://doi.org/10.1140/epjd/e2015-60464-1} {\bibfield  {journal} {\bibinfo
  {journal} {Eur. Phys. J. D}\ }\textbf {\bibinfo {volume} {69}},\ \bibinfo
  {pages} {279} (\bibinfo {year} {2015})}\BibitemShut {NoStop}%
\bibitem [{\citenamefont {Ac{\'{\i}}n}\ \emph {et~al.}(2018)\citenamefont
  {Ac{\'{\i}}n}, \citenamefont {Bloch}, \citenamefont {Buhrman}, \citenamefont
  {Calarco}, \citenamefont {Eichler}, \citenamefont {Eisert}, \citenamefont
  {Esteve}, \citenamefont {Gisin}, \citenamefont {Glaser}, \citenamefont
  {Jelezko}, \citenamefont {Kuhr}, \citenamefont {Lewenstein}, \citenamefont
  {Riedel}, \citenamefont {Schmidt}, \citenamefont {Thew}, \citenamefont
  {Wallraff}, \citenamefont {Walmsley},\ and\ \citenamefont
  {Wilhelm}}]{Acin2018}%
  \BibitemOpen
  \bibfield  {author} {\bibinfo {author} {\bibfnamefont {A.}~\bibnamefont
  {Ac{\'{\i}}n}}, \bibinfo {author} {\bibfnamefont {I.}~\bibnamefont {Bloch}},
  \bibinfo {author} {\bibfnamefont {H.}~\bibnamefont {Buhrman}}, \bibinfo
  {author} {\bibfnamefont {T.}~\bibnamefont {Calarco}}, \bibinfo {author}
  {\bibfnamefont {C.}~\bibnamefont {Eichler}}, \bibinfo {author} {\bibfnamefont
  {J.}~\bibnamefont {Eisert}}, \bibinfo {author} {\bibfnamefont
  {D.}~\bibnamefont {Esteve}}, \bibinfo {author} {\bibfnamefont
  {N.}~\bibnamefont {Gisin}}, \bibinfo {author} {\bibfnamefont {S.~J.}\
  \bibnamefont {Glaser}}, \bibinfo {author} {\bibfnamefont {F.}~\bibnamefont
  {Jelezko}}, \bibinfo {author} {\bibfnamefont {S.}~\bibnamefont {Kuhr}},
  \bibinfo {author} {\bibfnamefont {M.}~\bibnamefont {Lewenstein}}, \bibinfo
  {author} {\bibfnamefont {M.~F.}\ \bibnamefont {Riedel}}, \bibinfo {author}
  {\bibfnamefont {P.~O.}\ \bibnamefont {Schmidt}}, \bibinfo {author}
  {\bibfnamefont {R.}~\bibnamefont {Thew}}, \bibinfo {author} {\bibfnamefont
  {A.}~\bibnamefont {Wallraff}}, \bibinfo {author} {\bibfnamefont
  {I.}~\bibnamefont {Walmsley}},\ and\ \bibinfo {author} {\bibfnamefont
  {F.~K.}\ \bibnamefont {Wilhelm}},\ }\bibfield  {title} {\bibinfo {title} {The
  quantum technologies roadmap: a european community view},\ }\href
  {https://doi.org/10.1088/1367-2630/aad1ea} {\bibfield  {journal} {\bibinfo
  {journal} {New J. Phys.}\ }\textbf {\bibinfo {volume} {20}},\ \bibinfo
  {pages} {080201} (\bibinfo {year} {2018})}\BibitemShut {NoStop}%
\bibitem [{\citenamefont {Saywell}\ \emph {et~al.}(2018)\citenamefont
  {Saywell}, \citenamefont {Kuprov}, \citenamefont {Goodwin}, \citenamefont
  {Carey},\ and\ \citenamefont {Freegarde}}]{Saywell2018}%
  \BibitemOpen
  \bibfield  {author} {\bibinfo {author} {\bibfnamefont {J.~C.}\ \bibnamefont
  {Saywell}}, \bibinfo {author} {\bibfnamefont {I.}~\bibnamefont {Kuprov}},
  \bibinfo {author} {\bibfnamefont {D.}~\bibnamefont {Goodwin}}, \bibinfo
  {author} {\bibfnamefont {M.}~\bibnamefont {Carey}},\ and\ \bibinfo {author}
  {\bibfnamefont {T.}~\bibnamefont {Freegarde}},\ }\bibfield  {title} {\bibinfo
  {title} {Optimal control of mirror pulses for cold-atom interferometry},\
  }\href {https://doi.org/10.1103/physreva.98.023625} {\bibfield  {journal}
  {\bibinfo  {journal} {Phys. Rev. A}\ }\textbf {\bibinfo {volume} {98}},\
  \bibinfo {pages} {023625} (\bibinfo {year} {2018})}\BibitemShut {NoStop}%
\bibitem [{\citenamefont {Saywell}(2020)}]{SaywellThesis}%
  \BibitemOpen
  \bibfield  {author} {\bibinfo {author} {\bibfnamefont {J.~C.}\ \bibnamefont
  {Saywell}},\ }\emph {\bibinfo {title} {Optimal control of cold atoms for
  ultra-precise quantum sensors}},\ \href
  {http://eprints.soton.ac.uk/id/eprint/448869} {Ph.D. thesis},\ \bibinfo
  {school} {University of Southampton, UK} (\bibinfo {year} {2020})\BibitemShut
  {NoStop}%
\bibitem [{\citenamefont {Saywell}\ \emph {et~al.}(2020)\citenamefont
  {Saywell}, \citenamefont {Carey}, \citenamefont {Belal}, \citenamefont
  {Kuprov},\ and\ \citenamefont {Freegarde}}]{Saywell2020}%
  \BibitemOpen
  \bibfield  {author} {\bibinfo {author} {\bibfnamefont {J.}~\bibnamefont
  {Saywell}}, \bibinfo {author} {\bibfnamefont {M.}~\bibnamefont {Carey}},
  \bibinfo {author} {\bibfnamefont {M.}~\bibnamefont {Belal}}, \bibinfo
  {author} {\bibfnamefont {I.}~\bibnamefont {Kuprov}},\ and\ \bibinfo {author}
  {\bibfnamefont {T.}~\bibnamefont {Freegarde}},\ }\bibfield  {title} {\bibinfo
  {title} {Optimal control of raman pulse sequences for atom interferometry},\
  }\href {https://doi.org/10.1088/1361-6455/ab6df6} {\bibfield  {journal}
  {\bibinfo  {journal} {J. Phys. B: At. Mol. Opt. Phys.}\ }\textbf {\bibinfo
  {volume} {53}},\ \bibinfo {pages} {085006} (\bibinfo {year}
  {2020})}\BibitemShut {NoStop}%
\bibitem [{\citenamefont {Saywell}\ \emph {et~al.}(2021)\citenamefont
  {Saywell}, \citenamefont {Carey}, \citenamefont {Dedes}, \citenamefont
  {Kuprov},\ and\ \citenamefont {Freegarde}}]{Saywell2021}%
  \BibitemOpen
  \bibfield  {author} {\bibinfo {author} {\bibfnamefont {J.}~\bibnamefont
  {Saywell}}, \bibinfo {author} {\bibfnamefont {M.}~\bibnamefont {Carey}},
  \bibinfo {author} {\bibfnamefont {N.}~\bibnamefont {Dedes}}, \bibinfo
  {author} {\bibfnamefont {I.}~\bibnamefont {Kuprov}},\ and\ \bibinfo {author}
  {\bibfnamefont {T.}~\bibnamefont {Freegarde}},\ }\bibfield  {title} {\bibinfo
  {title} {Can optimised pulses improve the sensitivity of atom
  interferometers?},\ }in\ \href {https://doi.org/10.1117/12.2598991} {\emph
  {\bibinfo {booktitle} {Quantum Technology: Driving Commercialisation of an
  Enabling Science II}}},\ Vol.\ \bibinfo {volume} {11881},\ \bibinfo {editor}
  {edited by\ \bibinfo {editor} {\bibfnamefont {M.~J.}\ \bibnamefont
  {Padgett}}, \bibinfo {editor} {\bibfnamefont {K.}~\bibnamefont {Bongs}},
  \bibinfo {editor} {\bibfnamefont {A.}~\bibnamefont {Fedrizzi}},\ and\
  \bibinfo {editor} {\bibfnamefont {A.}~\bibnamefont {Politi}}},\ \bibinfo
  {organization} {International Society for Optics and Photonics}\ (\bibinfo
  {publisher} {SPIE},\ \bibinfo {year} {2021})\ pp.\ \bibinfo {pages}
  {83--92}\BibitemShut {NoStop}%
\bibitem [{\citenamefont {Zhang}\ \emph {et~al.}(2012)\citenamefont {Zhang},
  \citenamefont {Laflamme},\ and\ \citenamefont {Suter}}]{Zhang2012}%
  \BibitemOpen
  \bibfield  {author} {\bibinfo {author} {\bibfnamefont {J.}~\bibnamefont
  {Zhang}}, \bibinfo {author} {\bibfnamefont {R.}~\bibnamefont {Laflamme}},\
  and\ \bibinfo {author} {\bibfnamefont {D.}~\bibnamefont {Suter}},\ }\bibfield
   {title} {\bibinfo {title} {Experimental implementation of encoded logical
  qubit operations in a perfect quantum error correcting code},\ }\href
  {https://doi.org/10.1103/physrevlett.109.100503} {\bibfield  {journal}
  {\bibinfo  {journal} {Phys. Rev. Lett.}\ }\textbf {\bibinfo {volume} {109}},\
  \bibinfo {pages} {100503} (\bibinfo {year} {2012})}\BibitemShut {NoStop}%
\bibitem [{\citenamefont {Waldherr}\ \emph {et~al.}(2014)\citenamefont
  {Waldherr}, \citenamefont {Wang}, \citenamefont {Zaiser}, \citenamefont
  {Jamali}, \citenamefont {Schulte-Herbr{\"u}ggen}, \citenamefont {Abe},
  \citenamefont {Ohshima}, \citenamefont {Isoya}, \citenamefont {Du},
  \citenamefont {Neumann},\ and\ \citenamefont {Wrachtrup}}]{Waldherr2014}%
  \BibitemOpen
  \bibfield  {author} {\bibinfo {author} {\bibfnamefont {G.}~\bibnamefont
  {Waldherr}}, \bibinfo {author} {\bibfnamefont {Y.}~\bibnamefont {Wang}},
  \bibinfo {author} {\bibfnamefont {S.}~\bibnamefont {Zaiser}}, \bibinfo
  {author} {\bibfnamefont {M.}~\bibnamefont {Jamali}}, \bibinfo {author}
  {\bibfnamefont {T.}~\bibnamefont {Schulte-Herbr{\"u}ggen}}, \bibinfo {author}
  {\bibfnamefont {H.}~\bibnamefont {Abe}}, \bibinfo {author} {\bibfnamefont
  {T.}~\bibnamefont {Ohshima}}, \bibinfo {author} {\bibfnamefont
  {J.}~\bibnamefont {Isoya}}, \bibinfo {author} {\bibfnamefont {J.~F.}\
  \bibnamefont {Du}}, \bibinfo {author} {\bibfnamefont {P.}~\bibnamefont
  {Neumann}},\ and\ \bibinfo {author} {\bibfnamefont {J.}~\bibnamefont
  {Wrachtrup}},\ }\bibfield  {title} {\bibinfo {title} {Quantum error
  correction in a solid-state hybrid spin register},\ }\href
  {https://doi.org/10.1038/nature12919} {\bibfield  {journal} {\bibinfo
  {journal} {Nature}\ }\textbf {\bibinfo {volume} {506}},\ \bibinfo {pages}
  {204} (\bibinfo {year} {2014})}\BibitemShut {NoStop}%
\bibitem [{\citenamefont {Dolde}\ \emph {et~al.}(2014)\citenamefont {Dolde},
  \citenamefont {Bergholm}, \citenamefont {Wang}, \citenamefont {Jakobi},
  \citenamefont {Naydenov}, \citenamefont {Pezzagna}, \citenamefont {Meijer},
  \citenamefont {Jelezko}, \citenamefont {Neumann}, \citenamefont
  {Schulte-Herbr{\"{u}}ggen}, \citenamefont {Biamonte},\ and\ \citenamefont
  {Wrachtrup}}]{Dolde2014}%
  \BibitemOpen
  \bibfield  {author} {\bibinfo {author} {\bibfnamefont {F.}~\bibnamefont
  {Dolde}}, \bibinfo {author} {\bibfnamefont {V.}~\bibnamefont {Bergholm}},
  \bibinfo {author} {\bibfnamefont {Y.}~\bibnamefont {Wang}}, \bibinfo {author}
  {\bibfnamefont {I.}~\bibnamefont {Jakobi}}, \bibinfo {author} {\bibfnamefont
  {B.}~\bibnamefont {Naydenov}}, \bibinfo {author} {\bibfnamefont
  {S.}~\bibnamefont {Pezzagna}}, \bibinfo {author} {\bibfnamefont
  {J.}~\bibnamefont {Meijer}}, \bibinfo {author} {\bibfnamefont
  {F.}~\bibnamefont {Jelezko}}, \bibinfo {author} {\bibfnamefont
  {P.}~\bibnamefont {Neumann}}, \bibinfo {author} {\bibfnamefont
  {T.}~\bibnamefont {Schulte-Herbr{\"{u}}ggen}}, \bibinfo {author}
  {\bibfnamefont {J.}~\bibnamefont {Biamonte}},\ and\ \bibinfo {author}
  {\bibfnamefont {J.}~\bibnamefont {Wrachtrup}},\ }\bibfield  {title} {\bibinfo
  {title} {High-fidelity spin entanglement using optimal control},\ }\href
  {https://doi.org/10.1038/ncomms4371} {\bibfield  {journal} {\bibinfo
  {journal} {Nat. Commun.}\ }\textbf {\bibinfo {volume} {5}},\ \bibinfo {pages}
  {3371} (\bibinfo {year} {2014})}\BibitemShut {NoStop}%
\bibitem [{\citenamefont {Khaneja}\ \emph {et~al.}(2001)\citenamefont
  {Khaneja}, \citenamefont {Brockett},\ and\ \citenamefont
  {Glaser}}]{Khaneja2001}%
  \BibitemOpen
  \bibfield  {author} {\bibinfo {author} {\bibfnamefont {N.}~\bibnamefont
  {Khaneja}}, \bibinfo {author} {\bibfnamefont {R.}~\bibnamefont {Brockett}},\
  and\ \bibinfo {author} {\bibfnamefont {S.~J.}\ \bibnamefont {Glaser}},\
  }\bibfield  {title} {\bibinfo {title} {Time optimal control in spin
  systems},\ }\href {https://doi.org/10.1103/physreva.63.032308} {\bibfield
  {journal} {\bibinfo  {journal} {Phys. Rev. A}\ }\textbf {\bibinfo {volume}
  {63}},\ \bibinfo {pages} {032308} (\bibinfo {year} {2001})}\BibitemShut
  {NoStop}%
\bibitem [{\citenamefont {Skinner}\ \emph {et~al.}(2003)\citenamefont
  {Skinner}, \citenamefont {Reiss}, \citenamefont {Luy}, \citenamefont
  {Khaneja},\ and\ \citenamefont {Glaser}}]{Skinner2003}%
  \BibitemOpen
  \bibfield  {author} {\bibinfo {author} {\bibfnamefont {T.~E.}\ \bibnamefont
  {Skinner}}, \bibinfo {author} {\bibfnamefont {T.~O.}\ \bibnamefont {Reiss}},
  \bibinfo {author} {\bibfnamefont {B.}~\bibnamefont {Luy}}, \bibinfo {author}
  {\bibfnamefont {N.}~\bibnamefont {Khaneja}},\ and\ \bibinfo {author}
  {\bibfnamefont {S.~J.}\ \bibnamefont {Glaser}},\ }\bibfield  {title}
  {\bibinfo {title} {Application of optimal control theory to the design of
  broadband excitation pulses for high-resolution {NMR}},\ }\href
  {https://doi.org/10.1016/S1090-7807(03)00153-8} {\bibfield  {journal}
  {\bibinfo  {journal} {J. Magn. Reson.}\ }\textbf {\bibinfo {volume} {163}},\
  \bibinfo {pages} {8} (\bibinfo {year} {2003})}\BibitemShut {NoStop}%
\bibitem [{\citenamefont {Khaneja}\ \emph {et~al.}(2005)\citenamefont
  {Khaneja}, \citenamefont {Reiss}, \citenamefont {Kehlet}, \citenamefont
  {Schulte-Herbr{\"{u}}ggen},\ and\ \citenamefont {Glaser}}]{Khaneja2005}%
  \BibitemOpen
  \bibfield  {author} {\bibinfo {author} {\bibfnamefont {N.}~\bibnamefont
  {Khaneja}}, \bibinfo {author} {\bibfnamefont {T.}~\bibnamefont {Reiss}},
  \bibinfo {author} {\bibfnamefont {C.}~\bibnamefont {Kehlet}}, \bibinfo
  {author} {\bibfnamefont {T.}~\bibnamefont {Schulte-Herbr{\"{u}}ggen}},\ and\
  \bibinfo {author} {\bibfnamefont {S.~J.}\ \bibnamefont {Glaser}},\ }\bibfield
   {title} {\bibinfo {title} {Optimal control of coupled spin dynamics: design
  of {NMR} pulse sequences by gradient ascent algorithms},\ }\href
  {https://doi.org/10.1016/j.jmr.2004.11.004} {\bibfield  {journal} {\bibinfo
  {journal} {J. Magn. Reson.}\ }\textbf {\bibinfo {volume} {172}},\ \bibinfo
  {pages} {296} (\bibinfo {year} {2005})}\BibitemShut {NoStop}%
\bibitem [{\citenamefont {To{\v{s}}ner}\ \emph {et~al.}(2006)\citenamefont
  {To{\v{s}}ner}, \citenamefont {Glaser}, \citenamefont {Khaneja},\ and\
  \citenamefont {Nielsen}}]{Tosner2006}%
  \BibitemOpen
  \bibfield  {author} {\bibinfo {author} {\bibfnamefont {Z.}~\bibnamefont
  {To{\v{s}}ner}}, \bibinfo {author} {\bibfnamefont {S.~J.}\ \bibnamefont
  {Glaser}}, \bibinfo {author} {\bibfnamefont {N.}~\bibnamefont {Khaneja}},\
  and\ \bibinfo {author} {\bibfnamefont {N.~C.}\ \bibnamefont {Nielsen}},\
  }\bibfield  {title} {\bibinfo {title} {Effective {H}amiltonians by optimal
  control: {S}olid-state {NMR} double-quantum planar and isotropic dipolar
  recoupling},\ }\href {https://doi.org/10.1063/1.2366703} {\bibfield
  {journal} {\bibinfo  {journal} {J. Chem. Phys.}\ }\textbf {\bibinfo {volume}
  {125}},\ \bibinfo {pages} {184502} (\bibinfo {year} {2006})}\BibitemShut
  {NoStop}%
\bibitem [{\citenamefont {Vinding}\ \emph
  {et~al.}(2017{\natexlab{a}})\citenamefont {Vinding}, \citenamefont {Brenner},
  \citenamefont {Tse}, \citenamefont {Vellmer}, \citenamefont {Vosegaard},
  \citenamefont {Suter}, \citenamefont {St\"{o}cker},\ and\ \citenamefont
  {Maximov}}]{Vinding2017}%
  \BibitemOpen
  \bibfield  {author} {\bibinfo {author} {\bibfnamefont {M.~S.}\ \bibnamefont
  {Vinding}}, \bibinfo {author} {\bibfnamefont {D.}~\bibnamefont {Brenner}},
  \bibinfo {author} {\bibfnamefont {D.~H.~Y.}\ \bibnamefont {Tse}}, \bibinfo
  {author} {\bibfnamefont {S.}~\bibnamefont {Vellmer}}, \bibinfo {author}
  {\bibfnamefont {T.}~\bibnamefont {Vosegaard}}, \bibinfo {author}
  {\bibfnamefont {D.}~\bibnamefont {Suter}}, \bibinfo {author} {\bibfnamefont
  {T.}~\bibnamefont {St\"{o}cker}},\ and\ \bibinfo {author} {\bibfnamefont
  {I.~I.}\ \bibnamefont {Maximov}},\ }\bibfield  {title} {\bibinfo {title}
  {Application of the limited-memory quasi-{N}ewton algorithm for
  multi-dimensional, large flip-angle {RF} pulses at {7T}},\ }\href
  {https://doi.org/10.1007/s10334-016-0580-1} {\bibfield  {journal} {\bibinfo
  {journal} {Magn. Reson. Mater. Phy.}\ }\textbf {\bibinfo {volume} {30}},\
  \bibinfo {pages} {29} (\bibinfo {year} {2017}{\natexlab{a}})}\BibitemShut
  {NoStop}%
\bibitem [{\citenamefont {Reeth}\ \emph {et~al.}(2019)\citenamefont {Reeth},
  \citenamefont {Ratiney}, \citenamefont {Tse Ve~Koon}, \citenamefont {Tesch},
  \citenamefont {Grenier}, \citenamefont {Beuf}, \citenamefont {Glaser},\ and\
  \citenamefont {Sugny}}]{Reeth2019}%
  \BibitemOpen
  \bibfield  {author} {\bibinfo {author} {\bibfnamefont {E.~V.}\ \bibnamefont
  {Reeth}}, \bibinfo {author} {\bibfnamefont {H.}~\bibnamefont {Ratiney}},
  \bibinfo {author} {\bibfnamefont {K.}~\bibnamefont {Tse Ve~Koon}}, \bibinfo
  {author} {\bibfnamefont {M.}~\bibnamefont {Tesch}}, \bibinfo {author}
  {\bibfnamefont {D.}~\bibnamefont {Grenier}}, \bibinfo {author} {\bibfnamefont
  {O.}~\bibnamefont {Beuf}}, \bibinfo {author} {\bibfnamefont {S.~J.}\
  \bibnamefont {Glaser}},\ and\ \bibinfo {author} {\bibfnamefont
  {D.}~\bibnamefont {Sugny}},\ }\bibfield  {title} {\bibinfo {title} {A
  simplified framework to optimize {MRI} contrast preparation},\ }\href
  {https://doi.org/10.1002/mrm.27417} {\bibfield  {journal} {\bibinfo
  {journal} {Magn. Reson. Med.}\ }\textbf {\bibinfo {volume} {81}},\ \bibinfo
  {pages} {424} (\bibinfo {year} {2019})}\BibitemShut {NoStop}%
\bibitem [{\citenamefont {Soml{\'{o}}i}\ \emph {et~al.}(1993)\citenamefont
  {Soml{\'{o}}i}, \citenamefont {Kazakov},\ and\ \citenamefont
  {Tannor}}]{Somloi1993}%
  \BibitemOpen
  \bibfield  {author} {\bibinfo {author} {\bibfnamefont {J.}~\bibnamefont
  {Soml{\'{o}}i}}, \bibinfo {author} {\bibfnamefont {V.~A.}\ \bibnamefont
  {Kazakov}},\ and\ \bibinfo {author} {\bibfnamefont {D.~J.}\ \bibnamefont
  {Tannor}},\ }\bibfield  {title} {\bibinfo {title} {Controlled dissociation of
  {I}$_2$ via optical transitions between the {X} and {B} electronic states},\
  }\href {https://doi.org/10.1016/0301-0104(93)80108-l} {\bibfield  {journal}
  {\bibinfo  {journal} {Chem. Phys.}\ }\textbf {\bibinfo {volume} {172}},\
  \bibinfo {pages} {85} (\bibinfo {year} {1993})}\BibitemShut {NoStop}%
\bibitem [{\citenamefont {Zhu}\ \emph {et~al.}(1998)\citenamefont {Zhu},
  \citenamefont {Botina},\ and\ \citenamefont {Rabitz}}]{Zhu1998}%
  \BibitemOpen
  \bibfield  {author} {\bibinfo {author} {\bibfnamefont {W.}~\bibnamefont
  {Zhu}}, \bibinfo {author} {\bibfnamefont {J.}~\bibnamefont {Botina}},\ and\
  \bibinfo {author} {\bibfnamefont {H.}~\bibnamefont {Rabitz}},\ }\bibfield
  {title} {\bibinfo {title} {Rapidly convergent iteration methods for quantum
  optimal control of population},\ }\href {https://doi.org/10.1063/1.475576}
  {\bibfield  {journal} {\bibinfo  {journal} {J. Chem. Phys.}\ }\textbf
  {\bibinfo {volume} {108}},\ \bibinfo {pages} {1953} (\bibinfo {year}
  {1998})}\BibitemShut {NoStop}%
\bibitem [{\citenamefont {Maday}\ and\ \citenamefont
  {Turinici}(2003)}]{Maday2003}%
  \BibitemOpen
  \bibfield  {author} {\bibinfo {author} {\bibfnamefont {Y.}~\bibnamefont
  {Maday}}\ and\ \bibinfo {author} {\bibfnamefont {G.}~\bibnamefont
  {Turinici}},\ }\bibfield  {title} {\bibinfo {title} {New formulations of
  monotonically convergent quantum control algorithms},\ }\href
  {https://doi.org/10.1063/1.1564043} {\bibfield  {journal} {\bibinfo
  {journal} {J. Chem. Phys.}\ }\textbf {\bibinfo {volume} {118}},\ \bibinfo
  {pages} {8191–} (\bibinfo {year} {2003})}\BibitemShut {NoStop}%
\bibitem [{\citenamefont {Eitan}\ \emph {et~al.}(2011)\citenamefont {Eitan},
  \citenamefont {Mundt},\ and\ \citenamefont {Tannor}}]{Eitan2011}%
  \BibitemOpen
  \bibfield  {author} {\bibinfo {author} {\bibfnamefont {R.}~\bibnamefont
  {Eitan}}, \bibinfo {author} {\bibfnamefont {M.}~\bibnamefont {Mundt}},\ and\
  \bibinfo {author} {\bibfnamefont {D.~J.}\ \bibnamefont {Tannor}},\ }\bibfield
   {title} {\bibinfo {title} {Optimal control with accelerated convergence:
  {C}ombining the {K}rotov and quasi-{N}ewton methods},\ }\href
  {https://doi.org/10.1103/PhysRevA.83.053426} {\bibfield  {journal} {\bibinfo
  {journal} {Phys. Rev. A}\ }\textbf {\bibinfo {volume} {83}},\ \bibinfo
  {pages} {053426} (\bibinfo {year} {2011})}\BibitemShut {NoStop}%
\bibitem [{\citenamefont {Doria}\ \emph {et~al.}(2011)\citenamefont {Doria},
  \citenamefont {Calarco},\ and\ \citenamefont {Montangero}}]{Doria2011}%
  \BibitemOpen
  \bibfield  {author} {\bibinfo {author} {\bibfnamefont {P.}~\bibnamefont
  {Doria}}, \bibinfo {author} {\bibfnamefont {T.}~\bibnamefont {Calarco}},\
  and\ \bibinfo {author} {\bibfnamefont {S.}~\bibnamefont {Montangero}},\
  }\bibfield  {title} {\bibinfo {title} {Optimal control technique for
  many-body quantum dynamics},\ }\href
  {https://doi.org/10.1103/physrevlett.106.190501} {\bibfield  {journal}
  {\bibinfo  {journal} {Phys. Rev. Lett.}\ }\textbf {\bibinfo {volume} {106}},\
  \bibinfo {pages} {190501} (\bibinfo {year} {2011})}\BibitemShut {NoStop}%
\bibitem [{\citenamefont {Rach}\ \emph {et~al.}(2015)\citenamefont {Rach},
  \citenamefont {M{\"{u}}ller}, \citenamefont {Calarco},\ and\ \citenamefont
  {Montangero}}]{Rach2015}%
  \BibitemOpen
  \bibfield  {author} {\bibinfo {author} {\bibfnamefont {N.}~\bibnamefont
  {Rach}}, \bibinfo {author} {\bibfnamefont {M.~M.}\ \bibnamefont
  {M{\"{u}}ller}}, \bibinfo {author} {\bibfnamefont {T.}~\bibnamefont
  {Calarco}},\ and\ \bibinfo {author} {\bibfnamefont {S.}~\bibnamefont
  {Montangero}},\ }\bibfield  {title} {\bibinfo {title} {Dressing the
  chopped-random-basis optimization: {A} bandwidth-limited access to the
  trap-free landscape},\ }\href {https://doi.org/10.1103/physreva.92.062343}
  {\bibfield  {journal} {\bibinfo  {journal} {Phys. Rev. A}\ }\textbf {\bibinfo
  {volume} {92}},\ \bibinfo {pages} {062343} (\bibinfo {year}
  {2015})}\BibitemShut {NoStop}%
\bibitem [{\citenamefont {Ciaramella}\ \emph {et~al.}(2015)\citenamefont
  {Ciaramella}, \citenamefont {Borz{\`\i}}, \citenamefont {Dirr},\ and\
  \citenamefont {Wachsmuth}}]{Ciaramella2015}%
  \BibitemOpen
  \bibfield  {author} {\bibinfo {author} {\bibfnamefont {G.}~\bibnamefont
  {Ciaramella}}, \bibinfo {author} {\bibfnamefont {A.}~\bibnamefont
  {Borz{\`\i}}}, \bibinfo {author} {\bibfnamefont {G.}~\bibnamefont {Dirr}},\
  and\ \bibinfo {author} {\bibfnamefont {D.}~\bibnamefont {Wachsmuth}},\
  }\bibfield  {title} {\bibinfo {title} {Newton methods for the optimal control
  of closed quantum spin systems},\ }\href {https://doi.org/10.1137/140966988}
  {\bibfield  {journal} {\bibinfo  {journal} {{SIAM} J. Sci. Comput.}\ }\textbf
  {\bibinfo {volume} {37}},\ \bibinfo {pages} {A319} (\bibinfo {year}
  {2015})}\BibitemShut {NoStop}%
\bibitem [{\citenamefont {Machnes}\ \emph {et~al.}(2018)\citenamefont
  {Machnes}, \citenamefont {Ass{\'{e}}mat}, \citenamefont {Tannor},\ and\
  \citenamefont {Wilhelm}}]{Machnes2018}%
  \BibitemOpen
  \bibfield  {author} {\bibinfo {author} {\bibfnamefont {S.}~\bibnamefont
  {Machnes}}, \bibinfo {author} {\bibfnamefont {E.}~\bibnamefont
  {Ass{\'{e}}mat}}, \bibinfo {author} {\bibfnamefont {D.}~\bibnamefont
  {Tannor}},\ and\ \bibinfo {author} {\bibfnamefont {F.~K.}\ \bibnamefont
  {Wilhelm}},\ }\bibfield  {title} {\bibinfo {title} {Tunable, flexible, and
  efficient optimization of control pulses for practical qubits},\ }\href
  {https://doi.org/10.1103/physrevlett.120.150401} {\bibfield  {journal}
  {\bibinfo  {journal} {Phys. Rev. Lett.}\ }\textbf {\bibinfo {volume} {120}},\
  \bibinfo {pages} {150401} (\bibinfo {year} {2018})}\BibitemShut {NoStop}%
\bibitem [{\citenamefont {Lucarelli}(2018)}]{Lucarelli2018}%
  \BibitemOpen
  \bibfield  {author} {\bibinfo {author} {\bibfnamefont {D.}~\bibnamefont
  {Lucarelli}},\ }\bibfield  {title} {\bibinfo {title} {Quantum optimal control
  via gradient ascent in function space and the time-bandwidth quantum speed
  limit},\ }\href {https://doi.org/10.1103/PhysRevA.97.062346} {\bibfield
  {journal} {\bibinfo  {journal} {Phys. Rev. A}\ }\textbf {\bibinfo {volume}
  {97}},\ \bibinfo {pages} {062346} (\bibinfo {year} {2018})}\BibitemShut
  {NoStop}%
\bibitem [{\citenamefont {{Qui{\~n}ones Valles}}\ \emph
  {et~al.}(2019)\citenamefont {{Qui{\~n}ones Valles}}, \citenamefont {Dolgov},\
  and\ \citenamefont {Savostyanov}}]{QuinonesValles2019}%
  \BibitemOpen
  \bibfield  {author} {\bibinfo {author} {\bibfnamefont {D.}~\bibnamefont
  {{Qui{\~n}ones Valles}}}, \bibinfo {author} {\bibfnamefont {S.}~\bibnamefont
  {Dolgov}},\ and\ \bibinfo {author} {\bibfnamefont {D.}~\bibnamefont
  {Savostyanov}},\ }\bibfield  {title} {\bibinfo {title} {Tensor product
  approach to quantum control},\ }in\ \href
  {https://doi.org/10.1007/978-3-030-16077-7_29} {\emph {\bibinfo {booktitle}
  {Integral Methods in Science and Engineering: Analytic Treatment and
  Numerical Approximations}}},\ \bibinfo {editor} {edited by\ \bibinfo {editor}
  {\bibfnamefont {C.}~\bibnamefont {Constanda}}\ and\ \bibinfo {editor}
  {\bibfnamefont {P.}~\bibnamefont {Harris}}}\ (\bibinfo  {publisher}
  {Springer},\ \bibinfo {year} {2019})\ pp.\ \bibinfo {pages}
  {367--379}\BibitemShut {NoStop}%
\bibitem [{\citenamefont {Conolly}\ \emph {et~al.}(1986)\citenamefont
  {Conolly}, \citenamefont {Nishimura},\ and\ \citenamefont
  {A.}}]{Conolly1986}%
  \BibitemOpen
  \bibfield  {author} {\bibinfo {author} {\bibfnamefont {S.}~\bibnamefont
  {Conolly}}, \bibinfo {author} {\bibfnamefont {D.}~\bibnamefont {Nishimura}},\
  and\ \bibinfo {author} {\bibfnamefont {M.}~\bibnamefont {A.}},\ }\bibfield
  {title} {\bibinfo {title} {Optimal control solutions to the magnetic
  resonance selective excitation problem},\ }\href
  {https://doi.org/10.1109/TMI.1986.4307754} {\bibfield  {journal} {\bibinfo
  {journal} {IEEE Trans. Med. Imaging}\ }\textbf {\bibinfo {volume} {5}},\
  \bibinfo {pages} {106} (\bibinfo {year} {1986})}\BibitemShut {NoStop}%
\bibitem [{\citenamefont {Mao}\ \emph {et~al.}(1986)\citenamefont {Mao},
  \citenamefont {Mareci}, \citenamefont {Scott},\ and\ \citenamefont
  {Andrew}}]{Mao1986}%
  \BibitemOpen
  \bibfield  {author} {\bibinfo {author} {\bibfnamefont {J.}~\bibnamefont
  {Mao}}, \bibinfo {author} {\bibfnamefont {T.~H.}\ \bibnamefont {Mareci}},
  \bibinfo {author} {\bibfnamefont {K.~N.}\ \bibnamefont {Scott}},\ and\
  \bibinfo {author} {\bibfnamefont {E.~R.}\ \bibnamefont {Andrew}},\ }\bibfield
   {title} {\bibinfo {title} {Selective inversion radiofrequency pulses by
  optimal control},\ }\href {https://doi.org/10.1016/0022-2364(86)90016-8}
  {\bibfield  {journal} {\bibinfo  {journal} {J. Magn. Reson.}\ }\textbf
  {\bibinfo {volume} {70}},\ \bibinfo {pages} {310} (\bibinfo {year}
  {1986})}\BibitemShut {NoStop}%
\bibitem [{\citenamefont {Viola}\ and\ \citenamefont
  {Lloyd}(1998)}]{Viola1998}%
  \BibitemOpen
  \bibfield  {author} {\bibinfo {author} {\bibfnamefont {L.}~\bibnamefont
  {Viola}}\ and\ \bibinfo {author} {\bibfnamefont {S.}~\bibnamefont {Lloyd}},\
  }\bibfield  {title} {\bibinfo {title} {Dynamical suppression of decoherence
  in two-state quantum systems},\ }\href
  {https://doi.org/10.1103/PhysRevA.58.2733} {\bibfield  {journal} {\bibinfo
  {journal} {Phys. Rev. A}\ }\textbf {\bibinfo {volume} {58}},\ \bibinfo
  {pages} {2733} (\bibinfo {year} {1998})}\BibitemShut {NoStop}%
\bibitem [{\citenamefont {Viola}\ \emph {et~al.}(1999)\citenamefont {Viola},
  \citenamefont {Knill},\ and\ \citenamefont {Lloyd}}]{Viola1999}%
  \BibitemOpen
  \bibfield  {author} {\bibinfo {author} {\bibfnamefont {L.}~\bibnamefont
  {Viola}}, \bibinfo {author} {\bibfnamefont {E.}~\bibnamefont {Knill}},\ and\
  \bibinfo {author} {\bibfnamefont {S.}~\bibnamefont {Lloyd}},\ }\bibfield
  {title} {\bibinfo {title} {Dynamical decoupling of open quantum systems},\
  }\href {https://doi.org/10.1103/PhysRevLett.82.2417} {\bibfield  {journal}
  {\bibinfo  {journal} {Phys. Rev. Lett.}\ }\textbf {\bibinfo {volume} {82}},\
  \bibinfo {pages} {2417} (\bibinfo {year} {1999})}\BibitemShut {NoStop}%
\bibitem [{\citenamefont {de~Fouquieres}\ \emph {et~al.}(2011)\citenamefont
  {de~Fouquieres}, \citenamefont {Schirmer}, \citenamefont {Glaser},\ and\
  \citenamefont {Kuprov}}]{deFouquieres2011}%
  \BibitemOpen
  \bibfield  {author} {\bibinfo {author} {\bibfnamefont {P.}~\bibnamefont
  {de~Fouquieres}}, \bibinfo {author} {\bibfnamefont {S.~G.}\ \bibnamefont
  {Schirmer}}, \bibinfo {author} {\bibfnamefont {S.~J.}\ \bibnamefont
  {Glaser}},\ and\ \bibinfo {author} {\bibfnamefont {I.}~\bibnamefont
  {Kuprov}},\ }\bibfield  {title} {\bibinfo {title} {Second order gradient
  ascent pulse engineering},\ }\href
  {https://doi.org/10.1016/j.jmr.2011.07.023} {\bibfield  {journal} {\bibinfo
  {journal} {J. Magn. Reson.}\ }\textbf {\bibinfo {volume} {212}},\ \bibinfo
  {pages} {412} (\bibinfo {year} {2011})}\BibitemShut {NoStop}%
\bibitem [{\citenamefont {Goodwin}\ and\ \citenamefont
  {Kuprov}(2016)}]{Goodwin2016}%
  \BibitemOpen
  \bibfield  {author} {\bibinfo {author} {\bibfnamefont {D.~L.}\ \bibnamefont
  {Goodwin}}\ and\ \bibinfo {author} {\bibfnamefont {I.}~\bibnamefont
  {Kuprov}},\ }\bibfield  {title} {\bibinfo {title} {Modified
  {N}ewton-{R}aphson {GRAPE} methods for optimal control of spin systems},\
  }\href {https://doi.org/10.1063/1.4949534} {\bibfield  {journal} {\bibinfo
  {journal} {J. Chem. Phys.}\ }\textbf {\bibinfo {volume} {144}},\ \bibinfo
  {pages} {204107} (\bibinfo {year} {2016})}\BibitemShut {NoStop}%
\bibitem [{\citenamefont {Pechen}\ and\ \citenamefont
  {Tannor}(2011)}]{Pechen2011}%
  \BibitemOpen
  \bibfield  {author} {\bibinfo {author} {\bibfnamefont {A.~N.}\ \bibnamefont
  {Pechen}}\ and\ \bibinfo {author} {\bibfnamefont {D.~J.}\ \bibnamefont
  {Tannor}},\ }\bibfield  {title} {\bibinfo {title} {Are there traps in quantum
  control landscapes?},\ }\href
  {https://doi.org/10.1103/PhysRevLett.106.120402} {\bibfield  {journal}
  {\bibinfo  {journal} {Phys. Rev. Lett.}\ }\textbf {\bibinfo {volume} {106}},\
  \bibinfo {pages} {120402} (\bibinfo {year} {2011})}\BibitemShut {NoStop}%
\bibitem [{\citenamefont {Kobzar}\ \emph {et~al.}(2004)\citenamefont {Kobzar},
  \citenamefont {Skinner}, \citenamefont {Khaneja}, \citenamefont {Glaser},\
  and\ \citenamefont {Luy}}]{Kobzar2004}%
  \BibitemOpen
  \bibfield  {author} {\bibinfo {author} {\bibfnamefont {K.}~\bibnamefont
  {Kobzar}}, \bibinfo {author} {\bibfnamefont {T.~E.}\ \bibnamefont {Skinner}},
  \bibinfo {author} {\bibfnamefont {N.}~\bibnamefont {Khaneja}}, \bibinfo
  {author} {\bibfnamefont {S.~J.}\ \bibnamefont {Glaser}},\ and\ \bibinfo
  {author} {\bibfnamefont {B.}~\bibnamefont {Luy}},\ }\bibfield  {title}
  {\bibinfo {title} {Exploring the limits of broadband excitation and inversion
  pulses},\ }\href {https://doi.org/10.1016/j.jmr.2004.06.017} {\bibfield
  {journal} {\bibinfo  {journal} {J. Magn. Reson.}\ }\textbf {\bibinfo {volume}
  {170}},\ \bibinfo {pages} {236} (\bibinfo {year} {2004})}\BibitemShut
  {NoStop}%
\bibitem [{\citenamefont {Kobzar}\ \emph {et~al.}(2008)\citenamefont {Kobzar},
  \citenamefont {Skinner}, \citenamefont {Khaneja}, \citenamefont {Glaser},\
  and\ \citenamefont {Luy}}]{Kobzar2008}%
  \BibitemOpen
  \bibfield  {author} {\bibinfo {author} {\bibfnamefont {K.}~\bibnamefont
  {Kobzar}}, \bibinfo {author} {\bibfnamefont {T.~E.}\ \bibnamefont {Skinner}},
  \bibinfo {author} {\bibfnamefont {N.}~\bibnamefont {Khaneja}}, \bibinfo
  {author} {\bibfnamefont {S.~J.}\ \bibnamefont {Glaser}},\ and\ \bibinfo
  {author} {\bibfnamefont {B.}~\bibnamefont {Luy}},\ }\bibfield  {title}
  {\bibinfo {title} {Exploring the limits of broadband excitation and
  inversion: {II}. {Rf}-power optimized pulses},\ }\href
  {https://doi.org/10.1016/j.jmr.2008.05.023} {\bibfield  {journal} {\bibinfo
  {journal} {J. Magn. Reson.}\ }\textbf {\bibinfo {volume} {194}},\ \bibinfo
  {pages} {58} (\bibinfo {year} {2008})}\BibitemShut {NoStop}%
\bibitem [{\citenamefont {Goodwin}(2017)}]{GoodwinThesis}%
  \BibitemOpen
  \bibfield  {author} {\bibinfo {author} {\bibfnamefont {D.~L.}\ \bibnamefont
  {Goodwin}},\ }\emph {\bibinfo {title} {Advanced optimal control methods for
  spin systems}},\ \href {https://doi.org/10.5258/soton/t0003} {Ph.D. thesis},\
  \bibinfo  {school} {University of Southampton, UK} (\bibinfo {year}
  {2017})\BibitemShut {NoStop}%
\bibitem [{\citenamefont {To\u{s}ner}\ \emph {et~al.}(2021)\citenamefont
  {To\u{s}ner}, \citenamefont {Brandl}, \citenamefont {Blahut}, \citenamefont
  {Glaser},\ and\ \citenamefont {Reif}}]{Tosner2021}%
  \BibitemOpen
  \bibfield  {author} {\bibinfo {author} {\bibfnamefont {Z.}~\bibnamefont
  {To\u{s}ner}}, \bibinfo {author} {\bibfnamefont {M.~J.}\ \bibnamefont
  {Brandl}}, \bibinfo {author} {\bibfnamefont {J.}~\bibnamefont {Blahut}},
  \bibinfo {author} {\bibfnamefont {S.~J.}\ \bibnamefont {Glaser}},\ and\
  \bibinfo {author} {\bibfnamefont {B.}~\bibnamefont {Reif}},\ }\bibfield
  {title} {\bibinfo {title} {Maximizing efficiency of dipolar recoupling in
  solid-state {NMR} using optimal control sequences},\ }\href
  {https://doi.org/10.1126/sciadv.abj5913} {\bibfield  {journal} {\bibinfo
  {journal} {Sci. Adv.}\ }\textbf {\bibinfo {volume} {7}},\ \bibinfo {pages}
  {eabj5913} (\bibinfo {year} {2021})}\BibitemShut {NoStop}%
\bibitem [{\citenamefont {Kuprov}(2016)}]{Kuprov2016}%
  \BibitemOpen
  \bibfield  {author} {\bibinfo {author} {\bibfnamefont {I.}~\bibnamefont
  {Kuprov}},\ }\bibfield  {title} {\bibinfo {title} {{F}okker-{P}lanck
  formalism in magnetic resonance simulations},\ }\href
  {https://doi.org/10.1016/j.jmr.2016.07.005} {\bibfield  {journal} {\bibinfo
  {journal} {J. Magn. Reson.}\ }\textbf {\bibinfo {volume} {270}},\ \bibinfo
  {pages} {124} (\bibinfo {year} {2016})}\BibitemShut {NoStop}%
\bibitem [{\citenamefont {Vinding}\ \emph {et~al.}(2019)\citenamefont
  {Vinding}, \citenamefont {Skyum}, \citenamefont {Sangill},\ and\
  \citenamefont {Lund}}]{Vinding2019}%
  \BibitemOpen
  \bibfield  {author} {\bibinfo {author} {\bibfnamefont {M.~S.}\ \bibnamefont
  {Vinding}}, \bibinfo {author} {\bibfnamefont {B.}~\bibnamefont {Skyum}},
  \bibinfo {author} {\bibfnamefont {R.}~\bibnamefont {Sangill}},\ and\ \bibinfo
  {author} {\bibfnamefont {T.~E.}\ \bibnamefont {Lund}},\ }\bibfield  {title}
  {\bibinfo {title} {Ultrafast (milliseconds), multidimensional {RF} pulse
  design with deep learning},\ }\href {https://doi.org/10.1002/mrm.27740}
  {\bibfield  {journal} {\bibinfo  {journal} {Magn. Reson. Med.}\ }\textbf
  {\bibinfo {volume} {82}},\ \bibinfo {pages} {586} (\bibinfo {year}
  {2019})}\BibitemShut {NoStop}%
\bibitem [{\citenamefont {Vinding}\ \emph
  {et~al.}(2021{\natexlab{a}})\citenamefont {Vinding}, \citenamefont {Aigner},
  \citenamefont {Schmitter},\ and\ \citenamefont {Lund}}]{Vinding2021a}%
  \BibitemOpen
  \bibfield  {author} {\bibinfo {author} {\bibfnamefont {M.~S.}\ \bibnamefont
  {Vinding}}, \bibinfo {author} {\bibfnamefont {C.~S.}\ \bibnamefont {Aigner}},
  \bibinfo {author} {\bibfnamefont {S.}~\bibnamefont {Schmitter}},\ and\
  \bibinfo {author} {\bibfnamefont {T.~E.}\ \bibnamefont {Lund}},\ }\bibfield
  {title} {\bibinfo {title} {{DeepControl}: {2DRF} pulses facilitating
  inhomogeneity and {B0} off-resonance compensation in vivo at 7{T}},\ }\href
  {https://doi.org/10.1002/mrm.28667} {\bibfield  {journal} {\bibinfo
  {journal} {Magn. Reson. Med.}\ }\textbf {\bibinfo {volume} {85}},\ \bibinfo
  {pages} {3308} (\bibinfo {year} {2021}{\natexlab{a}})}\BibitemShut {NoStop}%
\bibitem [{\citenamefont {Vinding}\ \emph
  {et~al.}(2021{\natexlab{b}})\citenamefont {Vinding}, \citenamefont {Goodwin},
  \citenamefont {Kuprov},\ and\ \citenamefont {Lund}}]{Vinding2021b}%
  \BibitemOpen
  \bibfield  {author} {\bibinfo {author} {\bibfnamefont {M.~S.}\ \bibnamefont
  {Vinding}}, \bibinfo {author} {\bibfnamefont {D.~L.}\ \bibnamefont
  {Goodwin}}, \bibinfo {author} {\bibfnamefont {I.}~\bibnamefont {Kuprov}},\
  and\ \bibinfo {author} {\bibfnamefont {T.~E.}\ \bibnamefont {Lund}},\
  }\bibfield  {title} {\bibinfo {title} {Optimal control gradient precision
  trade-offs: {A}pplication to fast generation of {DeepControl} libraries for
  {MRI}},\ }\href {https://doi.org/10.1016/j.jmr.2021.107094} {\bibfield
  {journal} {\bibinfo  {journal} {J. Magn. Reson.}\ }\textbf {\bibinfo {volume}
  {333}},\ \bibinfo {pages} {107094} (\bibinfo {year}
  {2021}{\natexlab{b}})}\BibitemShut {NoStop}%
\bibitem [{\citenamefont {Goodwin}\ \emph {et~al.}(2020)\citenamefont
  {Goodwin}, \citenamefont {Koos},\ and\ \citenamefont {Luy}}]{Goodwin2020}%
  \BibitemOpen
  \bibfield  {author} {\bibinfo {author} {\bibfnamefont {D.~L.}\ \bibnamefont
  {Goodwin}}, \bibinfo {author} {\bibfnamefont {M.~R.~M.}\ \bibnamefont
  {Koos}},\ and\ \bibinfo {author} {\bibfnamefont {B.}~\bibnamefont {Luy}},\
  }\bibfield  {title} {\bibinfo {title} {Second order phase dispersion by
  optimised rotation pulses},\ }\href
  {https://doi.org/PhysRevResearch.2.033157} {\bibfield  {journal} {\bibinfo
  {journal} {Phys. Rev. Research}\ }\textbf {\bibinfo {volume} {2}},\ \bibinfo
  {pages} {033157} (\bibinfo {year} {2020})}\BibitemShut {NoStop}%
\bibitem [{\citenamefont {Haller}\ \emph {et~al.}(2022)\citenamefont
  {Haller}, \citenamefont {Goodwin},\ and\ \citenamefont {Luy}}]{Haller2022}%
  \BibitemOpen
  \bibfield  {author} {\bibinfo {author} {\bibfnamefont {J.~D.}\ \bibnamefont
  {Haller}}, \bibinfo {author} {\bibfnamefont {D.~L.}\ \bibnamefont
  {Goodwin}},\ and\ \bibinfo {author} {\bibfnamefont {B.}~\bibnamefont {Luy}},\
  }\bibfield  {title} {\bibinfo {title} {SORDOR pulses: expansion of the 
  B\"{o}hlen-Bodenhausen scheme for low-power broadband magnetic resonance},\ }\href
  {https://doi.org/10.5194/mr-3-53-2022} {\bibfield  {journal} {\bibinfo
  {journal} {Magn. Reson.}\ }\textbf {\bibinfo {volume} {3}},\ \bibinfo
  {pages} {53} (\bibinfo {year} {2022})}\BibitemShut {NoStop}%
\bibitem [{\citenamefont {Dempsey}\ \emph {et~al.}(2002)\citenamefont
  {Dempsey}, \citenamefont {Condon},\ and\ \citenamefont
  {Hadley}}]{Dempsey2002}%
  \BibitemOpen
  \bibfield  {author} {\bibinfo {author} {\bibfnamefont {M.~F.}\ \bibnamefont
  {Dempsey}}, \bibinfo {author} {\bibfnamefont {B.}~\bibnamefont {Condon}},\
  and\ \bibinfo {author} {\bibfnamefont {D.~M.}\ \bibnamefont {Hadley}},\
  }\bibfield  {title} {\bibinfo {title} {{MRI} safety review},\ }in\ \href
  {https://doi.org/10.1016/S0887-2171(02)90010-7} {\emph {\bibinfo {booktitle}
  {Seminars in Ultrasound, {CT} and {MRI}}}},\ Vol.~\bibinfo {volume} {23}\
  (\bibinfo {year} {2002})\ pp.\ \bibinfo {pages} {392--401}\BibitemShut
  {NoStop}%
\bibitem [{\citenamefont {Foroozandeh}\ and\ \citenamefont
  {Singh}(2021)}]{Foroozandeh2021}%
  \BibitemOpen
  \bibfield  {author} {\bibinfo {author} {\bibfnamefont {M.}~\bibnamefont
  {Foroozandeh}}\ and\ \bibinfo {author} {\bibfnamefont {P.}~\bibnamefont
  {Singh}},\ }\bibfield  {title} {\bibinfo {title} {Optimal control of spins by
  analytical {Lie} algebraic derivatives},\ }\href
  {https://doi.org/10.1016/j.automatica.2021.109611} {\bibfield  {journal}
  {\bibinfo  {journal} {Automatica}\ }\textbf {\bibinfo {volume} {129}},\
  \bibinfo {pages} {109611} (\bibinfo {year} {2021})}\BibitemShut {NoStop}%
\bibitem [{\citenamefont {Kallies}(2018)}]{KalliesThesis}%
  \BibitemOpen
  \bibfield  {author} {\bibinfo {author} {\bibfnamefont {W.}~\bibnamefont
  {Kallies}},\ }\emph {\bibinfo {title} {Concurrent optimization of robust
  refocused pulse sequences for magnetic resonance spectroscopy}},\ \href@noop
  {} {Ph.D. thesis},\ \bibinfo  {school} {Technische Universit\"{a}t
  M\"{u}nchen, DE} (\bibinfo {year} {2018})\BibitemShut {NoStop}%
\bibitem [{\citenamefont {Goodwin}\ and\ \citenamefont
  {Kuprov}(2015)}]{Goodwin2015}%
  \BibitemOpen
  \bibfield  {author} {\bibinfo {author} {\bibfnamefont {D.~L.}\ \bibnamefont
  {Goodwin}}\ and\ \bibinfo {author} {\bibfnamefont {I.}~\bibnamefont
  {Kuprov}},\ }\bibfield  {title} {\bibinfo {title} {Auxiliary matrix formalism
  for interaction representation transformations, optimal control, and spin
  relaxation theories},\ }\href {https://doi.org/10.1063/1.4928978} {\bibfield
  {journal} {\bibinfo  {journal} {J. Chem. Phys.}\ }\textbf {\bibinfo {volume}
  {143}},\ \bibinfo {pages} {084113} (\bibinfo {year} {2015})}\BibitemShut
  {NoStop}%
\bibitem [{\citenamefont {Siminovitch}(1997)}]{Siminovitch1997}%
  \BibitemOpen
  \bibfield  {author} {\bibinfo {author} {\bibfnamefont {D.~J.}\ \bibnamefont
  {Siminovitch}},\ }\bibfield  {title} {\bibinfo {title} {Rotations in {NMR}:
  Part {I}. {Euler-Rodrigues} parameters and quaternions},\ }\href
  {https://doi.org/10.1002/(SICI)1099-0534(1997)9:3<149::AID-CMR3>3.0.CO;2-Z}
  {\bibfield  {journal} {\bibinfo  {journal} {Concepts Magn. Reson.}\ }\textbf
  {\bibinfo {volume} {9}},\ \bibinfo {pages} {149} (\bibinfo {year}
  {1997})}\BibitemShut {NoStop}%
\bibitem [{\citenamefont {Glaser}\ \emph {et~al.}(1998)\citenamefont {Glaser},
  \citenamefont {Schulte-Herbr\"{u}ggen}, \citenamefont {Sieveking},
  \citenamefont {Schedletzky}, \citenamefont {Nielsen}, \citenamefont
  {S{\o}rensen},\ and\ \citenamefont {Griesinger}}]{Glaser1998}%
  \BibitemOpen
  \bibfield  {author} {\bibinfo {author} {\bibfnamefont {S.~J.}\ \bibnamefont
  {Glaser}}, \bibinfo {author} {\bibfnamefont {T.}~\bibnamefont
  {Schulte-Herbr\"{u}ggen}}, \bibinfo {author} {\bibfnamefont {M.}~\bibnamefont
  {Sieveking}}, \bibinfo {author} {\bibfnamefont {O.}~\bibnamefont
  {Schedletzky}}, \bibinfo {author} {\bibfnamefont {N.~C.}\ \bibnamefont
  {Nielsen}}, \bibinfo {author} {\bibfnamefont {O.~W.}\ \bibnamefont
  {S{\o}rensen}},\ and\ \bibinfo {author} {\bibfnamefont {C.}~\bibnamefont
  {Griesinger}},\ }\bibfield  {title} {\bibinfo {title} {Unitary control in
  quantum ensembles: Maximizing signal intensity in coherent spectroscopy},\
  }\href {https://doi.org/10.1126/science.280.5362.421} {\bibfield  {journal}
  {\bibinfo  {journal} {Science}\ }\textbf {\bibinfo {volume} {280}},\ \bibinfo
  {pages} {421} (\bibinfo {year} {1998})}\BibitemShut {NoStop}%
\bibitem [{\citenamefont {{Van Loan}}(1978)}]{VanLoan1978}%
  \BibitemOpen
  \bibfield  {author} {\bibinfo {author} {\bibfnamefont {C.~F.}\ \bibnamefont
  {{Van Loan}}},\ }\bibfield  {title} {\bibinfo {title} {Computing integrals
  involving the matrix exponential},\ }\href
  {https://doi.org/10.1109/tac.1978.1101743} {\bibfield  {journal} {\bibinfo
  {journal} {{IEEE} Trans. Automat. Contr.}\ }\textbf {\bibinfo {volume}
  {23}},\ \bibinfo {pages} {395} (\bibinfo {year} {1978})}\BibitemShut
  {NoStop}%
\bibitem [{\citenamefont {Goodwin}\ \emph {et~al.}(2022)\citenamefont
  {Goodwin}, \citenamefont {Singh},\ and\ \citenamefont
  {Foroozandeh}}]{Goodwin2022}%
  \BibitemOpen
  \bibfield  {author} {\bibinfo {author} {\bibfnamefont {D.~L.}\ \bibnamefont
  {Goodwin}}, \bibinfo {author} {\bibfnamefont {P.}~\bibnamefont {Singh}},\
  and\ \bibinfo {author} {\bibfnamefont {M.}~\bibnamefont {Foroozandeh}},\
  }\bibfield  {title} {\bibinfo {title} {Adaptive optimal control of entangled
  qubits}} (\bibinfo {year} {2022}),\ \bibinfo {note} {accepted for publication
  in Sci. Adv.}\BibitemShut {Stop}%
\bibitem [{\citenamefont {Kuprov}(2013)}]{Kuprov2013}%
  \BibitemOpen
  \bibfield  {author} {\bibinfo {author} {\bibfnamefont {I.}~\bibnamefont
  {Kuprov}},\ }\bibfield  {title} {\bibinfo {title} {Spin system trajectory
  analysis under optimal control pulses},\ }\href
  {https://doi.org/10.1016/j.jmr.2013.02.012} {\bibfield  {journal} {\bibinfo
  {journal} {J. Mag. Reson.}\ }\textbf {\bibinfo {volume} {233}},\ \bibinfo
  {pages} {107} (\bibinfo {year} {2013})}\BibitemShut {NoStop}%
\bibitem [{\citenamefont {Vinding}\ \emph
  {et~al.}(2017{\natexlab{b}})\citenamefont {Vinding}, \citenamefont {Gu\'{e}rin},
  \citenamefont {Vosegaard},\ and\ \citenamefont
  {Nielsen}}]{vinding_local_2017}%
  \BibitemOpen
  \bibfield  {author} {\bibinfo {author} {\bibfnamefont {M.}~\bibnamefont
  {Vinding}}, \bibinfo {author} {\bibfnamefont {B.}~\bibnamefont {Guérin}},
  \bibinfo {author} {\bibfnamefont {T.}~\bibnamefont {Vosegaard}},\ and\
  \bibinfo {author} {\bibfnamefont {N.}~\bibnamefont {Nielsen}},\ }\bibfield
  {title} {\bibinfo {title} {Local {SAR}, global {SAR}, and power-constrained
  large-flip-angle pulses with optimal control and virtual observation
  points},\ }\href {https://doi.org/10.1002/mrm.26086} {\bibfield  {journal}
  {\bibinfo  {journal} {Magnetic Resonance in Medicine}\ }\textbf {\bibinfo
  {volume} {77}},\ \bibinfo {pages} {374} (\bibinfo {year}
  {2017}{\natexlab{b}})}\BibitemShut {NoStop}%
\bibitem [{\citenamefont {Stockmann}\ and\ \citenamefont
  {Wald}(2018)}]{stockmann_vivo_2018}%
  \BibitemOpen
  \bibfield  {author} {\bibinfo {author} {\bibfnamefont {J.~P.}\ \bibnamefont
  {Stockmann}}\ and\ \bibinfo {author} {\bibfnamefont {L.~L.}\ \bibnamefont
  {Wald}},\ }\bibfield  {title} {\bibinfo {title} {In vivo {B$_0$} field
  shimming methods for {MRI} at 7 {T}},\ }\href
  {https://doi.org/10.1016/j.neuroimage.2017.06.013} {\bibfield  {journal}
  {\bibinfo  {journal} {NeuroImage}\ }\textbf {\bibinfo {volume} {168}},\
  \bibinfo {pages} {71} (\bibinfo {year} {2018})}\BibitemShut {NoStop}%
\end{thebibliography}
\end{document}